\documentclass[12pt]{article}

\newcommand{\real}{{\bf R}}
\newcommand{\complex}{{\bf C}}
\newcommand{\intplus}{{\bf N}}
\newcommand{\allint}{{\bf Z}}
\renewcommand{\div}{\mathop{\rm div}}

\renewcommand{\d}{\,{\rm d}}            
\newcommand{\D}{{\rm d}}                
\newcommand{\tr}{{\rm T}}

\newcommand{\cA}{{\cal A}}
\newcommand{\cB}{{\cal B}}
\newcommand{\cD}{{\cal D}}
\newcommand{\cL}{{\cal L}}

\newcommand{\cO}{{\cal O}}

\newcommand{\cS}{{\cal S}}

\newcommand{\cT}{{\cal T}}
\newcommand{\cK}{{\cal K}}

\newcommand{\uu}{{\bf u}}
\newcommand{\vv}{{\bf v}}
\newcommand{\xx}{{\bf x}}
\newcommand{\yy}{{\bf y}}

\newcommand{\eeta}{{\mbox{\boldmath$\eta$}}}
\newcommand{\xxi}{{\mbox{\boldmath$\xi$}}}

\newcommand{\one}{{\mbox{\boldmath 1}}}

\newtheorem{theorem}{Theorem}[section]
\newtheorem{lemma}[theorem]{Lemma}

\newtheorem{proposition}[theorem]{Proposition}
\newtheorem{corollary}[theorem]{Corollary}
\newtheorem{remark}[theorem]{Remark}

\newcommand{\reff}[1]{(\ref{#1})}

\newcommand{\inttwo}{\int_{\real^2}}
\newcommand{\proof}{{\noindent \bf Proof:\ }}
\newcommand{\half}{{\frac{1}{2}}}

\newcommand{\tf}{{\tilde{f}}}

\def\Re{{\rm Re\,}}

\def\ess{{\rm ess}}

\def\div{\mathop{\rm div}}

\def\epsilon{\varepsilon}
\def\phi{\varphi}

\def\build#1_#2^#3{\mathrel{
  \mathop{\kern 0pt#1}\limits_{#2}^{#3}}}
\def\QED{\mbox{}\hfill$\Box$}

\renewcommand{\:}{\thinspace :}

\newdimen\texpscorrection
\newdimen\figcenter
\def\figurewithtex #1 #2 #3 #4 #5\cr{\null
  {\goodbreak\figcenter=\hsize\relax
  \advance\figcenter by -#4truecm
  \divide\figcenter by 2
  \begin{figure}[hbt]
  \vskip #3truecm\noindent\hskip\figcenter
  \includegraphics{#1}{\hskip\texpscorrection\input #2 }
  \vskip 0.8truecm\noindent \vbox{\noindent {\footnotesize #5}}
  \end{figure}}}
\def\point#1 #2 #3 {\rlap{\kern #1 truecm
\raise #2 truecm \hbox{#3}}}

\usepackage{amsfonts}
\begin{document}

\title{Global stability of vortex solutions of the two-dimensional
Navier-Stokes equation}

\author{Thierry Gallay \\ Institut Fourier \\
Universit\'e de Grenoble I \\ BP 74 \\
38402 Saint-Martin d'H\`eres \\France
\and
C. Eugene Wayne \\
Department of Mathematics \\
and  Center for BioDynamics \\
Boston University \\
111 Cummington St.\\
Boston, MA 02215, USA}

\maketitle
\begin{abstract}
Both experimental and numerical studies of fluid motion indicate that
initially localized regions of vorticity tend to evolve into isolated
vortices and that these vortices then serve as organizing centers for
the flow. In this paper we prove that in two dimensions localized
regions of vorticity do evolve toward a vortex. More precisely we
prove that any solution of the two-dimensional Navier-Stokes equation
whose initial vorticity distribution is integrable converges to 
an explicit self-similar solution called ``Oseen's vortex''. 
This implies that the Oseen vortices are dynamically stable
for all values of the circulation Reynolds number, and our
approach also shows that these vortices are the only solutions
of the two-dimensional Navier-Stokes equation with a Dirac mass as 
initial vorticity. Finally, under slightly stronger assumptions on 
the vorticity distribution, we also give precise estimates on the 
rate of convergence toward the vortex.
\end{abstract}

\section{Introduction}\label{intro}

In this paper we consider the motion of an incompressible, viscous fluid in
two-dimensional Euclidean space.  The velocity of such a 
fluid is described by the Navier-Stokes equations
\begin{equation} \label{NS}
\frac{\partial \uu}{\partial t} +
  (\uu \cdot \nabla)\uu = \Delta \uu - \nabla p\ , \quad
  \nabla \cdot \uu = 0\ ,
\end{equation}
where $\uu = \uu(x,t) \in \real^2$ is the velocity field, $p = p(x,t) 
\in \real$ is the pressure field, and $x \in \real^2$, $t \ge 0$. 
For simplicity, the kinematic viscosity has been rescaled to $1$.

We prove two basic results about the solutions of \reff{NS}. First 
we show that for any initial velocity field whose vorticity is 
integrable, the solution of \reff{NS} with this initial velocity 
approaches an Oseen vortex, an explicit solution of \reff{NS} exhibited
below. As we also show, the Oseen vortex is in fact the 
unique solution of \reff{NS} with a Dirac mass as initial vorticity.
We then examine in more detail the approach toward the vortex by 
studying the spectrum of the linearized equation around the vortex 
solution. If we assume that the initial vorticity distribution lies 
in a weighted $L^2$ space, we can derive estimates on the 
spectrum of this linearized operator which allow us to prove
optimal bounds on the rate of convergence toward the vortex.

We now describe our results in more detail. As we have
argued in \cite{gallay:2001a} and \cite{gallay:2001b}
it is often easier to understand the asymptotics of solutions
of \reff{NS} by studying the evolution of the vorticity,
rather than the velocity.  This is especially true in 
two-dimensions where the vorticity is a scalar. One
can always then reconstruct the velocity field via the Biot-Savart
law. Taking the curl of \reff{NS} we find that the
vorticity $\omega = \partial_1 u_2 - \partial_2 u_1$ 
satisfies:
\begin{equation}\label{V}
  \frac{\partial \omega}{\partial t} + (\uu \cdot \nabla) \omega = 
  \Delta \omega\ , \quad x \in \real^2\ ,\quad t \ge 0\ .
\end{equation}
The velocity field $\uu$ is defined in terms of the vorticity 
via the Biot-Savart law
\begin{equation}\label{BS}
  \uu(x) = \frac{1}{2\pi} \inttwo 
  \frac{(\xx - \yy)^{\perp}}{|x -y|^2} \omega(y)\d y\ ,
  \quad x \in \real^2\ .
\end{equation}
Here and in the sequel, if $x=(x_1,x_2) \in \real^2$, we denote
$\xx = (x_1,x_2)^\tr$ and $\xx^{\perp} = (-x_2,x_1)^\tr$. 

The vorticity equation is globally well-posed in the space $L^1(\real^2)$.
In particular, the results of Ben-Artzi \cite{ben-artzi:1994},
Brezis \cite{brezis:1994} and Kato \cite{kato:1994} imply that:

\begin{theorem}\label{L1V2} For all initial data $\omega_0 \in L^1(\real^2)$, 
equation \reff{V} has a unique global solution $\omega \in 
C^0([0,\infty),L^1(\real^2)) \cap C^0((0,\infty),L^\infty(\real^2))$
such that $\omega(0) = \omega_0$. Moreover, for all $p \in [1,+\infty]$, 
there exists $C_p > 0$ such that
\begin{equation}\label{lpvorticity}
 |\omega(\cdot,t)|_p \le \frac{C_p |\omega_0|_1}{t^{1-\frac{1}{p}}}
 \ ,\quad t>0\ .
\end{equation}
\end{theorem}
\noindent Here and in the remainder of the paper $|\cdot|_p$ denotes 
the norm on $L^p(\real^2)$. If $\uu \in L^q(\real^2)^2$, we 
set $|u|_q = |\,|\uu|\,|_q$, where $|\uu| = (u_1^2+u_2^2)^{1/2}$.

\medskip
Among its other properties the semi-flow defined by \reff{V} in 
$L^1(\real^2)$ preserves mass, i.e.
\begin{equation}\label{masscons}
  \inttwo \omega(x,t) \d x = \inttwo \omega_0(x) \d x\ ,
  \quad t \ge 0\ .
\end{equation}
Furthermore, if the solution is sufficiently spatially localized
so that the first moments of the vorticity distribution are finite 
then these are also preserved:
\begin{equation}\label{momcons}
  \inttwo x_j \omega(x,t) \d x = \inttwo x_j \omega_0(x) \d x\ ,
  \quad t \ge 0\ , \quad j=1,2\ .
\end{equation}

It is important to realize that the solutions of \reff{V} given 
by Theorem~\ref{L1V2} correspond to {\it infinite energy} solutions
of the Navier-Stokes equations \reff{NS}. More precisely, if 
$\omega(x,t)$ is a solution of \reff{V} such that $\int \omega(x,t)
\d x \neq 0$, then the velocity field $\uu(x,t)$ given by 
\reff{BS} satisfies $|\uu(\cdot,t)|_2 = \infty$ for all $t$. 
Explicit examples of such infinite energy solutions are 
the so-called {\it Oseen vortices}: 
\begin{equation}\label{Oseen}
  \omega(x,t) = \frac{\alpha}{4 \pi t} \,e^{-|x|^2/(4t)}~,\quad
  \uu(x,t) = \frac{\alpha}{2\pi}\frac{\xx^\perp}{|x|^2}
  \Bigl(1 -  e^{-|x|^2/(4t)}\Bigr)~,
\end{equation}
where $|x|^2 = x_1^2 + x_2^2$ and $\alpha \in \real$ is a parameter
which is often referred to as the ``circulation Reynolds
number''. These solutions are ``trivial'' in the sense that
$\uu(x,t)\cdot\nabla\omega(x,t) \equiv 0$, so that \reff{V} reduces to
the linear heat equation.  However, they play a prominent role in the
long-time asymptotics of \reff{V}. Indeed, let
\begin{equation}\label{GvG}
  G(\xi) = \frac{1}{4\pi} \,e^{-|\xi|^2/4}~,\quad
  \vv^G(\xi) = \frac{1}{2\pi}\frac{\xxi^\perp}{|\xi|^2}
  \Bigl(1 -  e^{-|\xi|^2/4}\Bigr)~, \quad \xi \in \real^2\ .
\end{equation}
The following is the main result of this paper: 

\begin{theorem}\label{L1glob} If $\omega_0 \in L^1(\real^2)$,
 the solution $\omega(x,t)$ of \reff{V}
satisfies
\begin{equation}\label{L1conv}
  \lim_{t \to \infty} t^{1-\frac{1}{p}}\Big|\omega(\cdot,t) -
  \frac{\alpha}{t}G(\frac{\cdot}{\sqrt{t}})\Big|_p = 0~,
  ~~{\rm for}~~ 1 \le p \le \infty\ ,
\end{equation}
where $\alpha = \inttwo\omega_0(x)\d x$. 
If $\uu(x,t)$ is the solution of \reff{NS} obtained from
$\omega(x,t)$ via the Biot-Savart law \reff{BS}, then
\begin{equation}\label{L1conv2}
  \lim_{t \to \infty} t^{\half - \frac{1}{q}} \Big|
  \uu(\cdot,t) - \frac{\alpha}{\sqrt{t}}\vv^G(\frac{\cdot}{\sqrt{t}})
  \Big|_q = 0\ ,~~{\rm for}~~2 < q \le \infty \ .
\end{equation}
\end{theorem}

In other words, the solutions of \reff{V} in $L^1(\real^2)$ behave
asymptotically as the solutions of the linear heat equation
$\partial_t \omega = \Delta \omega$ with the same initial data.  For
small solutions, this result has been obtained by Giga and Kambe in
\cite{giga:1988b}, see also \cite{gallay:2001a}.  As was observed by
Carpio \cite{carpio:1994} (see also \cite{giga:2004}), there is a deep
connection between the asymptotics \reff{L1conv} and the uniqueness of
the fundamental solution of the vorticity equation. More precisely,
Carpio proved that \reff{L1conv} holds provided the Oseen vortex
\reff{Oseen} is the unique solution of \reff{V} with initial data
$\alpha\delta$, where $\delta$ is Dirac's measure. As is shown in
\cite{cottet:1986}, \cite{giga:1988}, \cite{kato:1994}, this is true
at least if $|\alpha|$ is sufficiently small. In this paper, we use a
different method which allows us to obtain \reff{L1conv} without any
restriction on $\alpha$. As a by-product of our analysis, we prove the
uniqueness of the solution of \reff{V} with a (large) Dirac mass as
initial condition:

\begin{proposition}\label{Unique}
Let $\omega \in C^0((0,T),L^1(\real^2)\cap L^\infty(\real^2))$
be a solution of \reff{V} which is bounded in $L^1(\real^2)$, 
and assume that $\omega(\cdot,t)$ (considered as a finite 
Radon measure on $\real^2$) converges weakly to $\alpha\delta$
as $t \to 0+$, for some $\alpha \in \real$. Then
$$
   \omega(x,t) = \frac{\alpha}{t} \,G\Bigl(\frac{x}{\sqrt{t}}
   \Bigr)~, \quad x \in \real^2~,\quad 0 < t < T~.
$$
\end{proposition}

\medskip
Theorem~\ref{L1glob} has a number of important consequences.
Recalling that $\alpha$ can be thought of as the Reynolds number of
the flow, Theorem~\ref{L1glob} says in more physical terms that the
Oseen vortices are globally stable for any value of this number.  In
contrast to many situations in hydrodynamics, such as the Poiseuille
or the Taylor-Couette flows, increasing the Reynolds number does not
produce any instability. From another point of view, our result is
compatible with the ``inverse cascade'' of energy in two-dimensional
turbulence theory.  In contrast to the situation in three dimensions
where energy injected into the system at large scales flows to smaller
and smaller scales until it is dissipated by viscosity, in two dimensions
both experimental and numerical results indicate that even for very
turbulent, high Reynolds number flows, there is a tendency for smaller
vortices to coalesce and form larger and larger vortices.  In this
context, Theorem~\ref{L1glob} says that in the whole space $\real^2$
this process continues until only a single vortex remains.

Another consequence of this result is that the Oseen vortices are the
only self-similar solutions of the Navier-Stokes equations in
$\real^2$ such that the vorticity field is integrable, see
\cite{rossi:2002} for a related result. For completeness, we mention
that these equations have many other self-similar solutions with
nonintegrable vorticities. Indeed, adapting the results of
\cite{cannone:1996} to the two-dimensional case, it is easy to verify
that the Cauchy problem for \reff{NS} is globally well-posed for small
data in a Besov space which contains homogeneous functions of degree
$-1$. For such initial data, the velocity $\uu(x,t)$ and the vorticity
$\omega(x,t)$ are automatically self-similar, due to scaling
invariance.  For instance, given any continuous map $\phi : S^1 \to
\real$, there exists $\epsilon > 0$ such that \reff{NS} has a
self-similar solution $\uu(x,t)$ with initial data $\uu_0(x) =
\epsilon x |x|^{-2} \phi(x/|x|)$. If $\phi$ is not constant, the
associated vorticity $\omega(x,t)$ decays like $1/|x|^2$ as $|x| \to
\infty$ (at least in some directions), so that $\omega(\cdot,t) \notin
L^1(\real^2)$.

\medskip
To prove Theorem~\ref{L1glob}, our strategy is to study 
a rescaled version of \reff{V} which is suggested by the
form of the vortex solution \reff{Oseen}. Thus, we introduce 
the ``scaling variables'' or ``similarity variables'':
$$
  \xi = \frac{x}{\sqrt{t}}\ , \qquad \tau = \log t\ .
$$
If $\omega(x,t)$ is a solution of \reff{V} and if $\uu(x,t)$ is the 
corresponding velocity field, we define new functions $w(\xi,\tau)$, 
$\vv(\xi,\tau)$ by
\begin{equation}\label{omega-w}
  \omega(x,t) = \frac{1}{t} \,w\Bigl(\frac{x}{\sqrt{t}},
  \log t\Bigr)\ , \quad 
  \uu(x,t) = \frac{1}{\sqrt{t}} \,\vv\Bigl(\frac{x}{\sqrt{t}},
  \log t\Bigr)\ .
\end{equation}
Then $w(\xi,\tau)$ satisfies the equation
\begin{equation}\label{SV}
  \partial_\tau  w + (\vv \cdot \nabla_\xi) w = \cL w\ ,
\end{equation}
where
\begin{equation}\label{Ldef}
  \cL w = \Delta_\xi w + \half (\xxi \cdot \nabla_\xi)w +w\ .
\end{equation}
The rescaled velocity $\vv$ is reconstructed from the rescaled 
vorticity $w$ by the Biot-Savart law:
\begin{equation}\label{SBS}
  \vv(\xi) = \frac{1}{2\pi} \inttwo
  \frac{(\xxi - \eeta)^{\perp}}{|\xi - \eta |^2} w(\eta) 
  \d\eta\ , \quad \xi \in \real^2\ .
\end{equation}
From \reff{omega-w} and Theorem~\ref{L1V2}, it is clear that
the Cauchy problem for \reff{SV} is globally well-posed in 
$L^1(\real^2)$. Remark that $w(\cdot,0) = \omega(\cdot,1)$, 
hence imposing initial conditions to \reff{SV} at time $\tau = 0$
corresponds to imposing initial conditions to \reff{V} at time 
$t=1$. This is of course harmless since \reff{V} is autonomous. 
Observe also that the Oseen vortices $\{\alpha G\}_{\alpha\in\real}$
are by construction a family of equilibria of \reff{SV}. 

To prove Theorem~\ref{L1glob}, we first study the long-time
asymptotics of solutions whose vorticity distribution is more strongly 
localized than is necessary just to be in $L^1$. For any $m \ge 0$, 
we introduce the weighted Lebesgue space $L^2(m)$ defined by
\begin{eqnarray*}
  L^2(m) &=& \bigl\{ f \in L^2(\real^2) \,|\, \|f\|_m < \infty \bigr\}\ ,
  \quad \hbox{where}\\ \nonumber
  \|f\|_m &=& \left(\inttwo(1+|\xi|^2)^m |f(\xi)|^2 \d\xi \right)^{1/2}
  = |b^m f|_2\ .
\end{eqnarray*}
Here and in the sequel, we denote $b(\xi) = (1+|\xi|^2)^{1/2}$. 
If $m > 1$, then $L^2(m) \hookrightarrow L^1(\real^2)$. In this case, 
we denote by $L^2_0(m)$ the closed subspace of $L^2(m)$ given by
\begin{equation}\label{L20mdef}
  L^2_0(m) = \Bigl\{w \in L^2(m)\,\Big|\, \inttwo w(\xi)\d\xi = 0
  \Bigr\}\ .
\end{equation}
For $\vv \in (L^2(m))^2$, we set $\| \vv \|_m = \|~ |\vv|~  \|_m$,
where $|\vv| = (v_1^2+v_2^2)^{1/2}$.

As we observed in \cite{gallay:2001a}, a crucial advantage of 
using the vorticity formulation of the Navier-Stokes equations
is that the spatial decay of the vorticity field is preserved 
under the evolution of \reff{V}. This remark is especially useful 
if one is interested in the long-time asymptotics 
since for parabolic equations there is a close relationship between 
the spatial and temporal decay of the solutions. The following
result shows that the Cauchy problem for \reff{SV} is globally
well-posed in the weighted space $L^2(m)$ if $m > 1$. 

\begin{theorem}[\cite{gallay:2001a}, Theorem 3.2]
\label{weightedL2} Suppose that $w_0 \in L^2(m)$ for some
$m > 1$.  Then \reff{SV} has a unique global solution $w \in 
C^0([0,\infty),L^2(m))$ with $w(0) = w_0$, and there exists
$C_1 = C_1(\|w_0\|_m) > 0$ such that
\begin{equation}\label{wglobound}
  \|w(\tau)\|_m \le C_1\ , \quad \tau \ge 0\ .
\end{equation}
Moreover, $C_1(\|w_0\|_m) \to 0$ as $\|w_0\|_m \to 0$. Finally, if 
$w_0 \in L^2_0(m)$, then $\inttwo w(\xi,\tau)\d\xi = 0$ for all 
$\tau \ge 0$, and $\displaystyle{\lim_{\tau\to\infty} \|w(\tau)\|_m = 0}$. 
\end{theorem}

Note that in contrast to this result the semi-flow defined by the
Navier-Stokes equation does not preserve the spatial localization of
the velocity field \cite{gallay:2001a}. For instance, if the initial
velocity $\uu_0$ lies in $L^2(m)$ for some $m > 2$, then in general
the solution $\uu(\cdot,t)$ of \reff{NS} with initial data $\uu_0$
will not be in $L^2(m)$ for $t > 0$. For a detailed study of the
localization properties of solutions of the Navier-Stokes and
vorticity equations, we refer to the recent work of Brandolese
\cite{brandolese:2001a}, \cite{brandolese:2001b}, \cite{brandolese:2003}.

\medskip
Our proof of Theorem~\ref{L1glob} begins with a proof that the
Oseen vortices attract all solutions of \reff{SV} with initial 
data in $L^2(m)$.

\begin{proposition}\label{global}
Let $m > 1$, $w_0 \in L^2(m)$, and let $w \in C^0([0,+\infty),L^2(m))$ 
be the solution of \reff{SV} with initial data $w_0$. Then 
$$
   \|w(\tau) - \alpha G\|_m \to 0 \quad \hbox{as }\tau \to +\infty~,
$$
where $\alpha = \inttwo w_0(\xi)\d\xi$.
\end{proposition}

Thus, any solution of the vorticity equation which is sufficiently
localized to be in $L^2(m)$ for $m>1$ will converge, as
time tends to infinity, toward one of the Oseen vortices,
regardless of how large the Reynolds number is. As is shown in 
Section~\ref{global_section}, this global convergence result can 
then be extended to all solutions in $L^1(\real^2)$. Returning 
to the original variables, we thus obtain Theorem~\ref{L1glob}
as a corollary.  

The proof of Proposition~\ref{global} is based on the existence 
of a pair of Lyapunov functionals for the rescaled vorticity 
equation \reff{SV}. The first Lyapunov function, which is just 
the $L^1$ norm, is nonincreasing due to the maximum principle.
It implies that the $\omega$-limit set of any solution must lie in 
the subset of solutions which are either everywhere positive or
everywhere negative. 

On the subset of positive solutions, we use a second Lyapunov
function which is motivated by a formal analogy between 
\reff{SV} and some kinetic models such as the 
Vlasov-Fokker-Planck equation. Given $w : \real^2 \to \real_+$, 
we set
$$ 
  H(w) \,=\, \inttwo w(\xi) \log \Bigl(\frac{w(\xi)}{G(\xi)}\Bigr) 
  \d\xi\ .
$$
This quantity is often called the {\it relative entropy} (or
relative information) of the vorticity distribution $w$ 
with respect to the Gaussian $G$, see e.g. \cite{villani:2002}. 
A direct calculation shows that $H$ is non-increasing along
the trajectories of \reff{SV}:
$$
  \frac{\D}{\D\tau}H(w) \,=\,
  -\inttwo w \Big|\nabla\log\Bigl(\frac{w}{G} 
  \Bigr)\Big|^2 \d\xi \le 0\ .
$$
More precisely, this formula shows that $H$ is strictly decreasing
except along the family of Oseen vortices, which play the role of the
Maxwellian distributions in kinetic theory. By LaSalle's invariance
principle, the $\omega$-limit set of any (nonnegative) solution of
\reff{SV} with total vorticity $\alpha$ is thus reduced to a single
point $\{\alpha G\}$, which proves Proposition~\ref{global}.

\medskip
We next investigate the rate at which the solution $w(\tau)$ of 
\reff{SV} approaches the vortex $\alpha G$ as $\tau \to \infty$.
This can be done by studying the linearization of \reff{SV}
at the vortex. Under the assumptions of Proposition~\ref{global}, 
there exists $\mu > 0$ such that $\|w(\tau) - \alpha G\|_m = 
\cO(e^{-\mu\tau})$ as $\tau \to \infty$. Convergence is 
thus exponential in the rescaled time $\tau = \log t$, 
hence algebraic in the original time $t$. As in the case
of the linear heat equation, the convergence rate satisfies 
$\mu < (m{-}1)/2$. This limitation originates from the essential 
spectrum of the linearized operator at the vortex, and is related
to the spatial decay of the vorticity. In addition, discrete
eigenvalues prevent convergence to be arbitrarily fast even for 
exponentially localized solutions. The following statement 
generalizes to arbitrary data in $L^2(m)$ the local results of 
\cite{gallay:2001a}. 

\begin{proposition}\label{local1}
Fix $m > 2$. For any $w_0 
\in L^2(m)$, the solution $w \in C^0( [0,+\infty),L^2(m))$ 
of \reff{SV} with initial data $w_0$ satisfies
\begin{equation}\label{loc1}
  \|w(\tau) - \alpha G \|_m = \cO(e^{-\tau/2})\ ,\quad 
  \hbox{as }\tau \to \infty\ ,
\end{equation}
where $\alpha = \inttwo w_0(\xi)\d\xi$. Moreover, if $m > 3$ 
and $\beta_1 = \beta_2 = 0$ where $\beta_i = \inttwo \xi_i w_0(\xi)
\d\xi$, then
\begin{equation}\label{loc2}
  \|w(\tau) - \alpha G \|_m = \cO(e^{-\tau})\ ,\quad 
  \hbox{as }\tau \to \infty\ .
\end{equation}
\end{proposition}

As we shall see in
Section~\ref{local_section}, if $m > 2$ and $(\beta_1,\beta_2) \neq
(0,0)$, then $\|w(\tau) - \alpha G \|_m$ decays exactly like
$e^{-\tau/2}$ as $\tau \to \infty$, so that \reff{loc1} is
sharp. Similarly, if $m > 3$, \reff{loc2} is sharp in the sense that
there is in general a correction to the Gaussian asymptotics decaying
exactly like $e^{-\tau}$ as $\tau \to \infty$. If
$\alpha \neq 0$, there is no loss of generality in assuming that
$\beta_1 = \beta_2 = 0$, since this can always been achieved by an
appropriate choice of the origin in the original variable $x \in
\real^2$ (see Section~\ref{local_section}). 
The situation is different if $\alpha = 0$, see
\cite{gallay:2001a}. In this case, if $(\beta_1,\beta_2) \neq (0,0)$,
the solution converges to zero at the rate $e^{-\tau/2}$, and the next
correction is $\cO(\tau e^{-\tau})$ due to secular terms in the
asymptotics.

If we reexpress estimate \reff{loc1} in terms of the original dependent
and independent variables, we find that
$$
  \Big|\omega(\cdot,t) - \frac{\alpha}{t}G(\frac{\cdot}{\sqrt{t}})
    \Big|_p = \cO(t^{-(\frac{3}{2}-\frac{1}{p})})\ , \quad
  \Big|\uu(\cdot,t) - \frac{\alpha}{\sqrt{t}}\vv^G(\frac{\cdot}{\sqrt{t}})
  \Big|_q = \cO(t^{-(1-\frac{1}{q})})\ , 
$$
as $t\to\infty$, which represents a considerable sharpening of the 
decay rates in \reff{L1conv}. Similarly, \reff{loc2} implies that 
the quantities above are $\cO(t^{-(2-\frac{1}{p})})$ and 
$\cO(t^{-(\frac{3}{2}-\frac{1}{q})})$, respectively.

The rest of this paper is organized as follows. In 
Section~\ref{compact_section}, we study the compactness 
properties of the solutions of \reff{SV} in both $L^2(m)$ 
and $L^1(\real^2)$. In Section~\ref{global_section}, 
we show that \reff{SV} has two Lyapunov functions, and
we prove Proposition~\ref{global}, Theorem~\ref{L1glob} and 
Proposition~\ref{Unique}. Finally, in Section~\ref{local_section}, 
we study the spectrum of the linearization of \reff{SV} at
the Oseen vortex and we obtain precise estimates of the rate
at which solutions of \reff{SV} approach the vortex. 
In particular, we prove Proposition~\ref{local1}.

\medskip\noindent{\bf Acknowledgements.} 
The first author is indebted to J. Dolbeault and, especially, 
to C. Villani for suggesting the beautiful idea of using the
Boltzmann entropy functional in the context of the 
two-dimensional Navier-Stokes equation. The research
of C.E.W. is supported in part by the NSF under grant 
DMS-0103915.


\section{Smoothing and compactness properties}\label{compact_section}

In this section we discuss some general properties of solutions of the
vorticity equation \reff{SV} which we will need to establish the
convergence results of the Section~\ref{global_section}.  To control
the nonlinear terms in \reff{V} or \reff{SV}, we will need estimates
on the velocity in terms of the vorticity. Let $\omega$ and $\uu$ be
related via the Biot-Savart law \reff{BS}. If $\omega \in
L^p(\real^2)$ for some $p \in (1,2)$, it follows from the classical
Hardy-Littlewood-Sobolev inequality that $\uu \in L^q(\real^2)^2$
where $\frac{1}{q} = \frac{1}{p} - \half$, and there exists $C > 0$
such that
\begin{equation}\label{HLS}
  |\uu|_q \le C |\omega|_p\ .
\end{equation}
Of course, a similar result holds if $w$ and $\vv$ are related 
via \reff{SBS}. Further estimates are collected in
(\cite{gallay:2001a}, Lemma~2.1 and Appendix B). 

\subsection{Compactness in $L^2(m)$}\label{compactL2}

If $w_0 \in L^2(m)$ for some $m > 1$, we know from Theorem~\ref{weightedL2}
that equation \reff{SV} has a unique global solution $w \in 
C^0([0,\infty),L^2(m))$ with initial data $w_0$. As is explained in 
\cite{gallay:2001a}, $w$ is in fact a solution of the associated 
integral equation
\begin{eqnarray}\label{SVint}
  w(\tau) &=& \cS(\tau) w_0 - \int_0^{\tau}\cS(\tau-s) \vv(s) \cdot \nabla 
    w(s) \d s \\ \nonumber
  &=& \cS(\tau) w_0 - \int_0^{\tau}e^{-\half(\tau-s)}\nabla \cdot
    \cS(\tau-s) \vv(s) w(s) \d s\ ,
\end{eqnarray}
where $\cS(\tau) = \exp(\tau\cL)$ is the $C_0$-semigroup generated by the 
operator $\cL$. Remark that, since $\cS(\tau)$ is not an analytic 
semigroup, the solution $w$ of \reff{SVint} is not (in general) 
a classical solution of \reff{SV}. In particular, $\tau \mapsto 
w(\tau)$ is not differentiable in $L^2(m)$. For later use, we recall
the following results (\cite{gallay:2001a}, Appendix A): 

\smallskip\noindent {\bf 1.} 
If $m > 1$, there exists $C > 0$ such that, for all $f \in L^2(m)$, 
\begin{equation}\label{L2mbdd}
  \|\cS(\tau)f\|_m \le C\|f\|_m\ , \quad   
  \|\nabla \cS(\tau)f\|_m \le \frac{C}{a(\tau)^{1/2}}\|f\|_m\ , \quad 
  \tau > 0\ ,
\end{equation}
where $a(\tau) = 1 - e^{-\tau}$. 

\smallskip\noindent {\bf 2.} 
If $0 < \mu \le 1/2$ and $m > 1 + 2\mu$, there exists $C > 0$ such that, 
for all $f \in L^2_0(m)$, 
\begin{equation}\label{L20mbdd}
  \|\cS(\tau)f\|_m \le C e^{-\mu\tau}\|f\|_m\ ,\quad
  \|\nabla \cS(\tau)f\|_m \le \frac{C e^{-\mu\tau}}{a(\tau)^{1/2}}
   \|f\|_m\ , \quad \tau > 0\ .
\end{equation}

\noindent {\bf 3.} Let $1 \le q \le p \le \infty$ and $T > 0$. For all $\alpha 
\in \intplus^2$ there exists $C > 0$ such that
\begin{equation}\label{LpLq}
  |b^m \partial^\alpha S(\tau)f|_p \le \frac{C}{a(\tau)^{
  (\frac{1}{q}-\frac{1}{p})+\frac{|\alpha|}{2}}} \,|b^m f|_q\ ,
  \quad 0 < \tau \le T\ ,
\end{equation}
where $b(\xi) = (1+|\xi|^2)^{1/2}$.

\medskip\noindent
From Theorem~\ref{weightedL2} we also know that there exists $K_1>0$
such that $\|w(\tau)\|_m \le K_1$ for all $\tau \ge 0$. The aim of this
section is to prove that the trajectory $\{w(\tau)\}_{\tau\ge 0}$ 
is in fact {\it relatively compact} in $L^2(m)$. This is because 
equation \reff{SV} is both regularizing and ``asymptotically 
confining'', in the sense that solutions of \reff{SV} in $L^2(m)$
are asymptotically (as $\tau \to +\infty$) well localized in 
space. The localization effect originates in the dilation 
term $\half \xi\cdot\nabla$ in the linear operator $\cL$, and
hence does not occur in the original vorticity equation \reff{V}.

We first prove that, for positive times, the solution $w(\tau)$ belongs 
to the weighted Sobolev space
\begin{equation}\label{H1def}
  H^1(m) = \bigl\{w \in L^2(m) \,|\, \partial_i w \in L^2(m) \hbox{ for }
   i = 1,2 \bigr\}\ ,\
\end{equation}
which we equip with the norm $\|w\|_{H^1(m)} = (\|w\|_m^2 + 
\|\nabla w\|_m^2)^{1/2}$.  

\begin{lemma}\label{SV_smooth} Let $w_0 \in L^2(m)$ with $m > 1$, and let 
$w \in C^0([0,\infty),L^2(m))$ be the solution of \reff{SV} with 
initial data $w_0$. Then there exists $K_2 > 0$ such that
\begin{equation}\label{smoothing_est}
  \|\nabla w(\tau)\|_m \le \frac{K_2}{a(\tau)^{1/2}}\ , \quad 
  \hbox{for all }\tau > 0\ ,
\end{equation} 
where $a(\tau)=1-e^{-\tau}$.
\end{lemma}

\proof Consider the Banach space $X = C^0([0,T],L^2(m)) \cap
C^0((0,T], H^1(m))$ equipped with the norm
$$
   \|w\|_X = \sup_{\tau \in [0,T]}\| w(\tau) \|_m + \sup_{\tau \in (0,T]}
   a(\tau)^{1/2}\|\nabla w(\tau)\|_m \ .
$$
We shall prove that there exist $T > 0$ and $K > 0$ such that, 
for all initial data $\tilde w_0 \in L^2(m)$ with $\|\tilde w_0\|_m 
\le K_1$, equation \reff{SVint} has a unique solution $\tilde w \in X$,
which satisfies $\|\tilde w\|_X \le K$. We then apply this result 
to $\tilde w_0 = w(nT/2)$ with $n \in \intplus$. By uniqueness, 
we have $\tilde w(\tau) = w(\tau + nT/2)$ for $\tau \in [0,T]$, 
hence
$$
   \sup_{0 < \tau \le T} a(\tau)^{1/2} \|\nabla w(\tau + nT/2)\|_m 
   \le K\ , \quad \hbox{for all }n \in \intplus\ ,
$$
which implies \reff{smoothing_est}. 

To prove that \reff{SVint} has a unique solution in $X$, 
we proceed as in (\cite{gallay:2001a}, Lemma~3.1). Given $w_1, w_2 
\in X$, we define
$$
   R(\tau) = \int_0^{\tau}\cS(\tau-s) \vv_1(s) \cdot \nabla w_2(s) \d s\ ,
  \quad 0 \le \tau \le T\ ,
$$
where $\vv_1(s)$ is the velocity field obtained from $w_1(s)$ 
via the Biot-Savart law \reff{SBS}. If $q \in (1,2)$, we know from 
\cite{gallay:2001a} that
$$
   \|R(\tau)\|_m \le C\Bigl(\int_0^\tau \frac{1}{a(s)^{1/q}}\d s\Bigr)
   \|w_1\|_X \|w_2\|_X\ .
$$
It remains to estimate
$$
  \nabla R(\tau) = \int_0^{\tau}\nabla \cS(\tau-s) \vv_1(s) \cdot 
  \nabla w_2(s) \d s\ .
$$
Applying \reff{LpLq} with $p = 2$ and $q$ as above, we obtain
$$
  \|\nabla \cS(\tau-s) \vv_1(s) \cdot \nabla w_2(s)\|_m
  \le \frac{C}{a(\tau-s)^{1/q}} |b^m \vv_1(s) \nabla w_2(s)|_q\ .
$$
Next, using H\"older's inequality, estimate \reff{HLS} and the fact 
that $L^2(m) \hookrightarrow L^q(\real^2)$ when $2/(m{+}1) < q \le 2$, 
we find
\begin{eqnarray}  \label{nonlinear_holder_bound}
  |b^m \vv_1(s) \cdot \nabla w_2(s)|_q \le
    |\vv_1(s)|_{\frac{2q}{2-q}} |b^m \nabla w_2(s)|_2
  &\le& C |w_1(s)|_q \|\nabla w_2(s)\|_m \\ \nonumber
  &\le& C \|w_1(s)\|_m \|\nabla w_2(s)\|_m\ .
\end{eqnarray}  
Thus
\begin{eqnarray} \nonumber
  a(\tau)^{1/2} \|\nabla R(\tau)\|_m &\le&
    \int_0^\tau \frac{C a(\tau)^{1/2}}{a(\tau-s)^{1/q}}
    \|w_1(s)\|_{m} \|\nabla w_2(s)\|_m \d s \\
  &\le& C \Bigl(\int_0^\tau \frac{a(\tau)^{1/2}}{a(\tau-s)^{1/q}a(s)^{1/2}}
    \d s\Bigr) \|w_1\|_X \|w_2\|_X  \ .
\end{eqnarray}
Summarizing, we have shown that $\|R\|_X \le C(T)\|w_1\|_X \|w_2\|_X$, 
where $C(T) \to 0$ as $T \to 0$. Using this bilinear estimate together
with the bound \reff{L2mbdd} on the linear semigroup, it is 
straightforward to verify by a fixed point argument that, if $T > 0$ 
is sufficiently small, then \reff{SVint} has a unique solution in 
$X$ with the desired properties. \QED

\medskip\noindent
Next, we decompose $w(\xi,\tau) = \alpha G(\xi) + R(\xi,\tau)$, 
where $G(\xi)$ is given by \reff{GvG} and $\alpha = \inttwo 
w(\xi,\tau)\d\xi$ is time-independent due to \reff{masscons}.
Then $R(\cdot,\tau) \in L^2_0(m)$ for all $\tau \ge 0$, where
$L^2_0(m)$ is defined in \reff{L20mdef}. Moreover, $R$ satisfies 
the evolution equation
\begin{equation}\label{R-evol}
  \partial_\tau R  \,=\, \cL R -\alpha \Lambda R  - N(R)\ ,
\end{equation}
where 
\begin{equation}\label{LamRdef}
  \Lambda R = \vv^G\cdot\nabla R + \vv^R\cdot\nabla G\ ,\quad
  N(R) = \vv^R\cdot\nabla R \ .
\end{equation}
Here $\vv^G$ is given by \reff{GvG} and $\vv^R$ denotes the velocity
field associated to the vorticity $R$ by the Biot-Savart law \reff{SBS}.
The corresponding integral equation is
\begin{equation}\label{Rint}
 R(\tau) \,=\, \cS(\tau) R_0 - \alpha \int_0^\tau \cS(\tau-s) 
  \Lambda R(s)  \d s - \int_0^\tau \cS(\tau-s) N(R(s))\d s \ .
\end{equation}
Equations \reff{SV} and \reff{R-evol} are clearly equivalent. 
In particular, Theorem~\ref{weightedL2} implies that, given
any $\alpha \in \real$, the Cauchy problem for \reff{R-evol} is 
globally well-posed in $L^2_0(m)$ if $m > 1$. 

\medskip
We now prove that positive trajectories of \reff{SV} in $L^2(m)$ are 
relatively compact.

\begin{lemma}\label{SV_compact} Let $w_0 \in L^2(m)$ with 
$m > 1$, and let $w \in C^0([0,\infty),L^2(m))$ be the solution 
of \reff{SV} with initial data $w_0$. Then the trajectory 
$\{w(\tau)\}_{\tau\ge 0}$ is relatively compact in $L^2(m)$. 
\end{lemma}

\proof Since $w \in C^0([0,\infty),L^2(m))$, it is sufficient 
to show that $\{w(\tau)\}_{\tau\ge 1}$ is relatively compact 
in $L^2(m)$. We decompose $w(\xi,\tau) = \alpha G(\xi) + R(\xi,\tau)$
as above, and consider the integral equation \reff{Rint} 
satisfied by $R$. 
Let $R_1(\tau) = R(\tau) - \cS(\tau)R_0$, so that $w(\tau) = \alpha G + 
\cS(\tau)R_0 + R_1(\tau)$. Since $R_0 \in L^2_0(m)$, it follows from
\reff{L20mbdd} that $\cS(\tau)R_0$ converges to zero in $H^1(m)$ 
as $\tau \to \infty$. In particular, since $\{w(\tau)\}_{\tau\ge 1}$
is bounded in $H^1(m)$ by Lemma~\ref{SV_smooth}, we see that 
$\{R_1(\tau)\}_{\tau\ge 1}$ is also bounded in $H^1(m)$. Now, 
we shall prove that $\{R_1(\tau)\}_{\tau\ge 0}$ is bounded in 
$L^2(m+1)$. By Rellich's criterion (see \cite{reed:1978}, 
Theorem~XIII.65), this implies that $\{R_1(\tau)\}_{\tau\ge 1}$
is relatively compact in $L^2(m)$. Since $\cS(\tau)R_0$ converges 
to zero as $\tau \to \infty$, it follows that $\{w(\tau)\}_{\tau\ge 1}$
is also relatively compact in $L^2(m)$, which is the desired result.

To prove the claim, we decompose $R_1(\tau) = R_2(\tau) + R_3(\tau) + 
R_4(\tau)$ where
\begin{eqnarray*}
  R_2(\tau) &=& -\alpha\int_0^\tau e^{-\half(\tau-s)}\nabla \cdot 
   \cS(\tau-s) \vv^G R(s) \d s\ ,\\
  R_3(\tau) &=& -\alpha\int_0^\tau e^{-\half(\tau-s)}\nabla \cdot 
   \cS(\tau-s) \vv^R(s) G \d s\ ,\\
  R_4(\tau) &=& -\int_0^\tau e^{-\half(\tau-s)}\nabla \cdot 
   \cS(\tau-s) \vv^R(s) R(s) \d s\ .
\end{eqnarray*}
To bound $R_2(\tau)$, we observe that $b\vv^G \in L^{\infty}(\real^2)$.
Using \reff{L2mbdd}, we thus find
\begin{eqnarray}\label{vG-R}\nonumber
 \|R_2(\tau)\|_{m+1} &=& |\alpha| \,\Big|b^{m+1} \int_0^{\tau} 
   e^{-\half(\tau-s)} \nabla \cdot \cS(\tau - s) \vv^G R(s) \Big|_2 \d s\\
 &\le& C |\alpha| \int_0^{\tau} e^{-\half (\tau-s)}
   \frac{1}{a(\tau-s)^{\half}} |b^{m+1} \vv^G R(s)|_2 \d s \\ \nonumber
 &\le& C |\alpha| \int_0^{\tau} e^{-\half (\tau-s)}
   \frac{1}{a(\tau-s)^{\half}} |b\vv^G|_{\infty} \|R(s)\|_m \d s
 \le C K_1^2\ ,
\end{eqnarray}
for all $\tau \ge 0$, where $K_1 = \sup_{\tau\ge 0}\|w(\tau)\|_m$. To
bound $R_3(\tau)$, we proceed in the same way. Using
H\"older's inequality, estimate \reff{HLS} and the embedding $L^2(m)
\hookrightarrow L^{4/3}$, we obtain
\begin{equation}\label{vR-G}
  |b^{m+1} \vv^R G|_2 \le |\vv^R|_4 |b^{m+1} G|_4 \le C |R|_{4/3} 
  \le C \|R\|_m\ ,
\end{equation}
so that $\|R_3(\tau)\|_{m+1} \le  C K_1^2$ for all $\tau \ge 0$.   
Finally, to bound $R_4(\tau)$ we fix $q \in (1,2)$ such that 
$q \ge 2/m$. Combining \reff{L2mbdd} and \reff{LpLq}, we obtain
\begin{equation}\label{vR-R}
  \|R_4(\tau)\|_{m+1} \le C \int_0^{\tau} e^{-\half (\tau-s)}
  \frac{1}{a(\tau-s)^{1/q}} |b^{m+1} \vv^R(s) R(s) |_q \d s\ .
\end{equation}
By H\"older's inequality, $|b\vv^R b^m R|_q \le |b \vv^R|_{\frac{2q}{2-q}}
\|R\|_m$. Applying Proposition~B.1(2) in \cite{gallay:2001a}, 
we also obtain $|b \vv^R|_{\frac{2q}{2-q}} \le C \|R\|_{\frac{2}{q}} \le 
C \|R\|_m$. Inserting these bounds into \reff{vR-R}, we find that 
$\|R_4(\tau)\|_{m+1} \le C K_1^2$ for all $\tau \ge 0$. Summarizing, 
we have shown that $\|R_1(\tau)\|_{m+1} \le C K_1^2$ 
for all $\tau \ge 0$, which concludes the proof. \QED 

\medskip
Finally, we show that negative or complete trajectories of 
\reff{SV} in $L^2(m)$ that are bounded in $L^2(m)$ for some
$m > 1$ are also relatively compact. 

\begin{lemma}\label{SV_neg} Assume that $m > 1$ and that 
$w \in C^0(\real,L^2(m))$ is a solution of \reff{SV} which is 
bounded in $L^2(m)$. Then $\{w(\tau)\}_{\tau\in\real}$ is 
relatively compact in $L^2(m)$.
\end{lemma}

\proof By assumption, there exists $K_1 > 0$ such that 
$\|w(\tau)\|_m \le K_1$ for all $\tau \in \real$. As in the 
proof of Lemma~\ref{SV_compact}, we decompose $w(\xi,\tau) = 
\alpha G(\xi) + R(\xi,\tau)$, where $\alpha = \inttwo w(\xi,\tau)\d\xi$.
The remainder $R(\tau)$ satisfies the integral equation
\begin{equation}\label{Rinteg}
  R(\tau) \,=\, \cS(\tau-\tau_0) R(\tau_0) - \alpha \int_{\tau_0}^\tau 
  \cS(\tau-s) \Lambda R(s) \d s - \int_{\tau_0}^\tau \cS(\tau-s) 
  N(R(s))\d s \ ,
\end{equation}
for all $\tau_0 < \tau$. Since $R(\tau_0) \in L^2_0(m)$ and 
$\|R(\tau_0)\|_m \le C K_1$ for all $\tau_0 \in \real$, it 
follows from \reff{L20mbdd} that $\|\cS(\tau-\tau_0)R(\tau_0)\|_m
\to 0$ as $\tau_0 \to -\infty$. Moreover, proceeding as in the
proof of Lemma~\ref{SV_compact} and using the analogues of 
estimates \reff{vG-R}, \reff{vR-R}, it is easy to see that 
both integrals in the right-hand side of \reff{Rinteg} have
a limit in $L^2(m)$ (and even in $L^2(m+1)$) as $\tau_0 \to
-\infty$. Thus, we have the representation
\begin{equation}\label{Rinteg2}
  R(\tau) \,=\, - \alpha \int_{-\infty}^\tau 
  \cS(\tau-s) \Lambda R(s) \d s - \int_{-\infty}^\tau \cS(\tau-s) 
  N(R(s))\d s \ , \quad \tau \in \real\ ,
\end{equation}
which implies that $\|R(\tau)\|_{m+1} \le C K_1^2$ for all $\tau \in \real$.
This shows that $\{w(\tau)\}_{\tau\in\real}$ is bounded in $L^2(m+1)$.
On the other hand, it follows from Lemma~\ref{SV_smooth} that 
$\{w(\tau)\}_{\tau\in\real}$ is bounded in $H^1(m)$, hence 
$\{w(\tau)\}_{\tau\in\real}$ is relatively compact in $L^2(m)$
by Rellich's criterion. 
\QED

\begin{remark}\label{bootstrap} By a bootstrap argument, it is 
clear from the proof of Lemma~\ref{SV_neg} that $\{w(\tau)\}_{\tau\in\real}$
is bounded in $H^k(m')$ for all $k \in \intplus$ and all $m' \in
\intplus$. In other words, the trajectory $\{w(\tau)\}_{\tau\in\real}$
is bounded in the Schwartz space $\cS(\real^2)$. 
\end{remark}

\subsection{Compactness in $L^1(\real^2)$}\label{compactL1}

We now study the compactness properties of the solutions of 
\reff{SV} in $L^1(\real^2)$. We first recall two important 
estimates for the solutions of the original vorticity equation
\reff{V}, and use the change of variables \reff{omega-w} 
to obtain the corresponding bounds on the solutions of \reff{SV}. 

The first is a smoothing estimate, see for instance 
\cite{kato:1994}. Under the assumptions of Theorem~\ref{L1V2}, 
the vorticity $\omega(x,t)$ satisfies for all $p \in (1,\infty]$
\begin{equation}\label{lpgrad}
 |\nabla \omega(\cdot,t)|_p \le \frac{C_p'}{t^{\frac{3}{2}
  -\frac{1}{p}}}\ ,\quad t>0\ ,
\end{equation}
where $C_p'$ depends only on $|\omega_0|_1$. 
The second is a nice pointwise estimate due to Carlen and Loss, see
Theorem~3 in \cite{carlen:1996}. For any $\beta \in (0,1)$, there 
exists $C_\beta > 0$ (depending on $|\omega_0|_1$) such that
\begin{equation}\label{pointbdd}
  |\omega(x,t)| \le C_\beta \inttwo \frac{1}{t}
  \exp\Bigl(-\beta\frac{|x-y|^2}{4t}\Bigr)|\omega_0(y)|\d y~,
  \quad x \in \real^2~, \quad t > 0\ . 
\end{equation}

Assume now that $w_0 \in L^1(\real^2)$, and let $w \in C^0([0,\infty),
L^1(\real^2))$ be the solution of \reff{SV} with initial data $w_0$.
Then the function $\omega(x,t)$ defined by \reff{omega-w} is 
a solution of \reff{V} on the time interval $[1,\infty)$ 
with initial data $\omega(x,1) = w_0(x)$. Applying \reff{lpgrad}
and returning to the rescaled variables, we obtain
\begin{equation}\label{lpgrad2}
 |\nabla w(\cdot,\tau)|_p \le \frac{C_p'}{a(\tau)^{\frac{3}{2}
  -\frac{1}{p}}}\ ,\quad \tau > 0\ ,
\end{equation}
where $a(\tau) = 1 - e^{-\tau}$. Similarly, we deduce from \reff{pointbdd}
that
\begin{equation}\label{CL1}
  |w(\xi,\tau)| \le C_\beta \inttwo \frac{1}{a(\tau)}\,
  \exp\Bigl(-\beta\frac{|\xi-y e^{-\tau/2}|^2}{4a(\tau)}
  \Bigr)|w_0(y)|\d y~, \quad \xi \in \real^2~, \quad \tau > 0\ .
\end{equation}
With these estimates at hand, we are now ready to prove the 
main result of this section. 

\begin{lemma}\label{L1comp}
Let $w_0 \in L^1(\real^2)$, and let $w \in C^0([0,\infty),L^1(\real^2))$
be the solution of \reff{SV} with initial data $w_0$. Then 
$\{w(\tau)\}_{\tau \ge 0}$ is relatively compact in $L^1(\real^2)$.
\end{lemma}

\proof
Again, it is sufficient to show that $\{w(\tau)\}_{\tau \ge 1}$
is relatively compact. Fix $R > 0$. From \reff{CL1} one has 
$|w(\xi,\tau)| \le w_1(\xi,\tau) + w_2(\xi,\tau)$ where
\begin{eqnarray*}
  w_1(\xi,\tau) &=& C_\beta \int_{|y|\le R} \frac{1}{a(\tau)}\,
    \exp\Bigl(-\beta\frac{|\xi-y e^{-\tau/2}|^2}{4a(\tau)}
    \Bigr)|w_0(y)|\d y~, \\
  w_2(\xi,\tau) &=& C_\beta \int_{|y|\ge R} \frac{1}{a(\tau)}\,
    \exp\Bigl(-\beta\frac{|\xi-y e^{-\tau/2}|^2}{4a(\tau)}
    \Bigr)|w_0(y)|\d y~. \\
\end{eqnarray*}
If $|\xi| \ge 2R$, then $|\xi-y e^{-\tau/2}| \ge |\xi| -|y| \ge 
|\xi|/2$ whenever $|y| \le R$ and $\tau \ge 0$. It follows that
$$
   w_1(\xi,\tau) \le C_\beta |w_0|_1 \,\frac{1}{a(\tau)}\,
    \exp\Bigl(-\beta\frac{|\xi|^2}{16a(\tau)}\Bigr)\ , 
   \quad |\xi| \ge 2R\ ,
$$
hence
$$
  \int_{|\xi|\ge 2R} w_1(\xi,\tau) \d\xi \,\le\, C_\beta |w_0|_1
  \int_{|z|\ge 2R} \exp\Bigl(-\beta\frac{|z|^2}{16}\Bigr)\d z
  \,=:\, \epsilon_1(R)\ . 
$$
Moreover, using Fubini's theorem, we find
$$
  \inttwo w_2(\xi,\tau) \d\xi \,\le\, C_\beta 
  \frac{4\pi}{\beta} \int_{|y|\ge R} |w_0(\xi)| \d\xi 
  \,=:\, \epsilon_2(R)\ . 
$$
Thus we have shown
\begin{equation}\label{CL2}
  \sup_{\tau > 0}\int_{|\xi|\ge 2R}|w(\xi,\tau)|\d\xi 
  \,\le\, \epsilon_1(R) + \epsilon_2(R) \to 0 \quad \hbox{as }
  R \to \infty~.
\end{equation}

On the other hand, using \reff{lpgrad2} with $p = \infty$, we
see that
$$
  K := \sup_{\tau\ge 1} \sup_{\xi \in \real^2} 
  |\nabla w(\xi,\tau)| < \infty~.
$$
Fix $\epsilon > 0$. According to \reff{CL2}, there exists 
$R > 1$ such that 
$$
  \sup_{\tau \ge 1}\int_{|\xi|\ge R-1}|w(\xi,\tau)|\d\xi 
  \le \epsilon/3\ .
$$
Let $\delta = \min(1,\epsilon(3K\pi R^2)^{-1})$. If $y \in 
\real^2$ satisfies $|y| \le \delta$, then for all $\tau \ge 1$
$$
  \int_{|\xi|\ge R} |w(\xi-y,\tau) - w(\xi,\tau)| \d\xi \,\le\, 
  2 \int_{|\xi|\ge R-1} |w(\xi,\tau)| \d\xi \le \frac{2\epsilon}{3}\ ,
$$
and
$$
  \int_{|\xi|\le R} |w(\xi-y,\tau) - w(\xi,\tau)| \d\xi \,\le\, 
  \pi R^2 |y| \sup_{|\xi|\le R+1}|\nabla w(\xi,\tau)| \le
  \pi R^2 \delta K \le \frac{\epsilon}{3}\ ,
$$
so that $\inttwo |w(\xi-y,\tau) - w(\xi,\tau)| \d\xi \le \epsilon$. 
Thus we have shown
\begin{equation}\label{CL3}
  \sup_{\tau\ge 1} \sup_{|y| \le \delta} \inttwo
  |w(\xi-y,\tau) - w(\xi,\tau)| \d\xi \to 0 
  \quad \hbox{as }\delta \to 0~.
\end{equation}
By the Riesz criterion (\cite{reed:1978}, Theorem~XIII.66), 
it follows from \reff{CL2}, \reff{CL3} that $\{w(\tau)\}_{\tau\ge 1}$
is relatively compact in $L^1(\real^2)$, which is the desired result. 
\QED

\medskip
To conclude this subsection, we remark that Lemma~\ref{SV_neg} has no
analogue in $L^1(\real^2)$, namely negative trajectories of \reff{SV}
that are bounded in $L^1(\real^2)$ need not be relatively compact in
that space. To see this, let $\omega_0 \in L^1(\real^2)$ and let
$\omega(x,t)$ be the solution of \reff{V} given by
Theorem~\ref{L1V2}. If $w(\xi,\tau) = e^\tau \omega(\xi
e^{\tau/2},e^\tau)$, then $\{w(\cdot,\tau)\}_{\tau\le 0}$ is a
bounded negative trajectory in $L^1(\real^2)$.  However,
$w(\cdot,\tau)$ is ``evanescent'' as $\tau \to -\infty$, hence
$\{w(\cdot,\tau)\}_{\tau\le 0}$ is not relatively compact in
$L^1(\real^2)$ unless $\omega_0 \equiv 0$.  

\subsection{Preservation of positivity}\label{positivity}

A more qualitative property that will be essential for our analysis is 
the fact that solutions of \reff{V}, and hence of \reff{SV}, satisfy a 
maximum principle. We state this property in the original 
variables, and for a generalized version of \reff{V} where 
the velocity and vorticity field are not necessarily connected 
to each other. This generalization will be useful in the 
next section. 

\begin{proposition}\label{max_prin} Assume that 
$\tilde{\uu} \in C^0_b(\real^2 \times [0,\infty),\real^2)$
and that $\omega \in C^2_b(\real^2 \times [0,\infty),\real)$
is a solution of
\begin{equation}
  \partial_t \omega(x,t) + \tilde{\uu}(x,t) \cdot \nabla 
  \omega(x,t) = \Delta \omega(x,t)\ , \quad x \in \real^2\,,~
  t \ge 0\ .
\end{equation}
If $\omega(x,0) \ge 0$ for all $x \in \real^2$, then either 
$\omega(x,t) \equiv 0$ or $\omega(x,t) > 0$ for all 
$x \in \real^2$ and all $t > 0$. 
\end{proposition}

\proof This classical result is obtained for instance by combining 
Theorems~3.5 and 3.10 in the book of Protter and Weinberger 
\cite{protter:1967}.
\QED

\medskip As a corollary, we obtain that \reff{V} preserves 
positivity. The same property holds for \reff{SV} and is 
proved using the change of variables \reff{omega-w}. 

\begin{corollary}\label{Smax-prin}
Assume that $\omega_0 \in L^1(\real^2)$ satisfies $\omega_0(x) 
\ge 0$ for almost all $x \in \real^2$, and that $\omega_0(x)$
does not vanish almost everywhere. Then the solution of 
\reff{V} given by Theorem~\ref{L1V2} satisfies $\omega(x,t) > 0$
for all $x \in \real^2$ and all $t > 0$. 
\end{corollary}

\proof
If $\omega_0 \in \cS(\real^2)$ and $\omega_0(x) \ge 0$ 
for all $x \in \real^2$, the solution $\omega(x,t)$ of 
\reff{V} and the corresponding velocity field $\uu(x,t)$ 
satisfy the assumptions of Proposition~\ref{max_prin}, 
hence $\omega(x,t) \ge 0$ for all $x \in \real^2$ and all $t > 0$. 
Since solutions of \reff{V} depend continuously on the 
initial data in $L^1(\real^2)$, the same result holds for 
any solution of \reff{V} with initial data $\omega_0 \in 
L^1(\real^2)$ such that $\omega_0(x) \ge 0$ almost everywhere.
Moreover, given any $t_0 > 0$, this solution satisfies 
$\omega \in C^2_b(\real^2 \times [t_0,\infty),\real)$ and
the corresponding velocity field satisfies $\uu \in C^0_b(\real^2 
\times [t_0,\infty),\real^2)$. Applying Proposition~\ref{max_prin} 
again, we deduce that $\omega(x,t) > 0$ for all $x\in \real^2$
and all $t > t_0$, unless $\omega(x,t_0) \equiv 0$. Since this
is true for any $t_0 > 0$, we conclude that $\omega(x,t) > 0$ for 
all $x\in \real^2$ and all $t > 0$, unless $\omega_0(x)$ 
vanishes almost everywhere. 
\QED


\section{Global convergence results}\label{global_section}

This section is devoted to the proofs of Proposition~\ref{global}, 
Theorem~\ref{L1glob} and Proposition~\ref{Unique}. The argument 
relies on the compactness properties of the previous section, 
and uses the crucial fact that system \reff{SV} has two Lyapunov 
functions which we introduce now. 

\subsection{A pair of Lyapunov functions}\label{Lyapunov}

Let $\Phi : L^1(\real^2) \to \real_+$ be the continuous 
function defined by
\begin{equation} \label{Phi-def}
  \Phi(w) \,=\, \inttwo |w(\xi)| \d\xi\ ,
\end{equation}
and let 
$$
   \Sigma \,=\, \left\{w \in L^1(\real^2) \,\Big|\, \inttwo
   |w(\xi)|\d \xi = \Big| \inttwo w(\xi)\d \xi\Big|\right\}~.
$$
In words, a function $w \in L^1(\real^2)$ belongs to $\Sigma$ if 
and only if $w(\xi)$ has (almost everywhere) a constant sign. 
Remark that, by Corollary~\ref{Smax-prin}, the set $\Sigma$ is 
positively invariant under the evolution of \reff{SV}. 

We first show that $\Phi$ is a Lyapunov function for the
semiflow of \reff{SV}. More precisely, $\Phi$ is strictly 
decreasing along the trajectories of \reff{SV} except on the 
invariant set $\Sigma$ where $\Phi$ is constant. 

\begin{lemma}\label{PhiLyap}
Let $w_0 \in L^1(\real^2)$, and let $w \in C^0([0,\infty),L^1(\real^2))$
be the solution of \reff{SV} with initial data $w_0$. Then 
$\Phi(w(\tau)) \le \Phi(w_0)$ for all $\tau \ge 0$. Moreover 
$\Phi(w(\tau)) = \Phi(w_0)$ for all $\tau \ge 0$ if and only if 
$w_0 \in \Sigma$. 
\end{lemma}

\proof If $w_0 \in \Sigma$, then $w(\tau) \in \Sigma$ for all
$\tau \ge 0$. Using \reff{masscons}, we thus find
$$
   \Phi(w(\tau)) \,=\, \Big|\inttwo w(\xi,\tau)\d \xi\Big| \,=\,
   \Big|\inttwo w_0(\xi)\d \xi\Big| = \Phi(w_0)\ ,
   \quad \hbox{for all } \tau \ge 0\ .
$$
Assume now that $w_0 \notin \Sigma$. Then $w_0 = w_0^+ - w_0^-$, 
where
$$
   w_0^+(\xi) = \max(w_0(\xi),0) \ge 0\ , \quad 
   w_0^-(\xi) = -\min(w_0(\xi),0) \ge 0\ . 
$$
By assumption, both $w_0^+$ and $w_0^-$ are nonzero on a set of 
positive Lebesgue measure. Let $w_1$ and $w_2$ be solutions of 
\begin{eqnarray}\label{w1-w2} \nonumber
  \partial_{\tau} w_1 + \vv \cdot \nabla w_1 &=& \cL w_1\ , 
    \quad \tau \ge 0\ ,  \\
  \partial_{\tau} w_2 + \vv \cdot \nabla w_2 &=& \cL w_2\ ,
    \quad \tau \ge 0\ ,  
\end{eqnarray}
with initial data $w_1(0) = w_0^+$, $w_2(0) = w_0^-$, where
$\vv(\xi,\tau)$ is the velocity field associated with the 
solution $w(\xi,\tau)$ of \reff{SV}. Following the proof of
Theorem~\ref{L1V2} (see \cite{kato:1994}), one verifies that
$w_1, w_2 \in C^0([0,\infty),L^1(\real^2))$ and that
\begin{equation}\label{massw1w2}
   \inttwo w_1(\xi,\tau)\d\xi \,=\, \inttwo w_0^+(\xi)\d\xi\ , \quad
   \inttwo w_2(\xi,\tau)\d\xi \,=\, \inttwo w_0^-(\xi)\d\xi\ , 
\end{equation}
for all $\tau \ge 0$. Moreover, using Proposition~\ref{max_prin}
and a density argument as in the proof of Corollary~\ref{Smax-prin}, 
it is straightforward to verify that $w_1(\xi,\tau) > 0$ and 
$w_2(\xi,\tau) > 0$ for all $\xi \in \real^2$ and all $\tau > 0$. 
Now, by construction we have $w(\xi,\tau) = w_1(\xi,\tau) -
w_2(\xi,\tau)$, hence
$$
  |w(\xi,\tau)| \,=\, |w_1(\xi,\tau) - w_2(\xi,\tau)| 
  \,<\, w_1(\xi,\tau) + w_2(\xi,\tau)\ ,
$$
for all $\xi \in \real^2$ and all $\tau > 0$. Integrating over
$\real^2$ and using \reff{massw1w2}, we obtain
\begin{eqnarray*}
  \inttwo |w(\xi,\tau)|\d\xi &<& \inttwo (w_1(\xi,\tau) + 
    w_2(\xi,\tau)) \d\xi \\
  &=& \inttwo (w_0^+(\xi) + w_0^-(\xi)) \d\xi \,=\, 
    \inttwo |w_0(\xi)|\d\xi~, \quad \tau > 0\ . 
\end{eqnarray*}
This shows that $\Phi(w(\tau)) < \Phi(w_0)$ for all $\tau > 0$. 
\QED

\medskip
Next we fix $m > 3$ and we consider solutions of \reff{SV} in 
the invariant cone $L^2(m) \cap \Sigma_+$, where
$$
   \Sigma_+ \,=\, \{w \in L^1(\real^2) \,|\, w(\xi) \ge 0
   \hbox{ almost everywhere} \}\ .
$$
We define $H : L^2(m) \cap \Sigma_+ \to \real$ by
\begin{equation}\label{Hdef} 
  H(w) \,=\, \inttwo w(\xi) \log \Bigl(\frac{w(\xi)}{G(\xi)}\Bigr) 
  \d\xi\ .
\end{equation}
Since
\begin{eqnarray*}
 w\log \Bigl(\frac{w}{G}\Bigr) &=& 
  \left(\frac{w}{G} \log\Bigl(\frac{w}{G} \Bigr)\right) G
  \,\ge\, -\frac{1}{e}\, G\ ,\\ 
 w\log \Bigl(\frac{w}{G}\Bigr) &=& w \log(4\pi w) + \frac{|\xi|^2}{4} 
 \, w \,\le\, C w^2 + \frac{|\xi|^2}{4}\, w\ ,
\end{eqnarray*}
it is clear that $H$ is well-defined and bounded from below by
$-1/e$. Moreover, using for instance the inequality
$$
   |w_1\log w_1 - w_2\log w_2| \,\le\, C\Bigl(|w_1-w_2|^{1/2}
   + |w_1-w_2| (w_1^{1/2} + w_2^{1/2})\Bigr)\ ,
$$
one verifies that $H$ is continuous on $L^2(m) \cap \Sigma_+$ 
(equipped with the topology of $L^2(m)$). We now show that $H$ is 
indeed a Lyapunov function for the semiflow defined by \reff{SV} 
on $L^2(m) \cap \Sigma_+$. 

\begin{lemma}\label{HLyap} Assume that $w_0 \in L^2(m) \cap \Sigma_+$ 
with $m > 3$, and let $w \in C^0([0,\infty),L^2(m))$
be the solution of \reff{SV} with initial data $w_0$. Then 
$H(w(\tau)) \le H(w_0)$ for all $\tau \ge 0$. Moreover 
$H(w(\tau)) = H(w_0)$ for all $\tau \ge 0$ if and only if 
$w_0 = \alpha G$ for some $\alpha \ge 0$.
\end{lemma}

\proof If $w_0 = \alpha G$ for some $\alpha \ge 0$, then 
$w(\tau) = \alpha G$ for all $\tau \ge 0$, hence obviously 
$H(w(\tau)) = H(w_0)$ for all $\tau \ge 0$. Assume now that $w_0$ 
is not a multiple of $G$ (in particular, $w_0 \neq 0$). 
Then the solution $w(\xi,\tau)$ of \reff{SV} is smooth and strictly
positive for all $\tau > 0$. We claim that $\tau \mapsto H(w(\tau))$
is differentiable for $\tau > 0$, and that
\begin{equation}\label{H-dot}
  \frac{\D}{\D\tau}H(w(\tau)) \,=\, -I(w(\tau))\ , \quad \tau > 0\ ,
\end{equation}
where
\begin{equation}\label{Idef}
  I(w) = \inttwo w(\xi) \Big|\nabla\log\Bigl(\frac{w(\xi)}{G(\xi)} 
  \Bigr)\Big|^2 \d\xi \ge 0\ .
\end{equation}
Remark that $I(w)$ vanishes if and only if $w$ is proportional
to $G$. Thus, under the assumptions above, it is clear that 
$I(w(\tau)) > 0$ at least for $\tau > 0$ sufficiently small, 
hence $H(w(\tau)) < H(w_0)$ for all $\tau > 0$. 

Thus, all that remains is to prove \reff{H-dot}. Assume first
that $w_0$ belongs to the Schwartz space $\cS(\real^2)$. 
Then $w \in C^1([0,\infty),\cS(\real^2))$ is a classical 
solution of \reff{SV} in $\cS(\real^2)$. Moreover, $w(\xi,\tau)$
satisfies a {\em Gaussian lower bound} for any $\tau > 0$, 
see \cite{osada:1987} or (\cite{giga:1988}, Theorem~3.1). 
More precisely, there exist positive constants $\gamma$ and 
$C_\gamma$ (depending only on $|w_0|_1$) such that, for all
$\xi \in \real^2$ and all $\tau> 0$,
\begin{eqnarray}\nonumber
  w(\xi,\tau) &\ge& \frac{C_\gamma}{a(\tau)} \inttwo 
   \exp\Bigl(-\gamma\frac{|\xi-y e^{-\tau/2}|^2}{2a(\tau)}
   \Bigr) w_0(y) \d y~, \\ \label{lowerbd}
  &\ge& \frac{C_\gamma}{a(\tau)} \,\exp\Bigl(-\gamma
   \frac{|\xi|^2}{a(\tau)}\Bigr)\inttwo \exp\Bigl(-\gamma 
   \frac{|y|^2}{a(\tau)}\Bigr) w_0(y) \d y~,
\end{eqnarray}
see also \reff{CL1}. Using these properties, it is straightforward to 
verify that $\tau \mapsto H(w(\tau))$ is differentiable for 
$\tau > 0$, and that
$$
  \frac{\D}{\D\tau}H(w(\tau)) \,=\,
  \inttwo \Bigl(1+\log\frac{w}{G}\Bigr) \partial_\tau 
  w\d\xi = \inttwo \Bigl(1+\log\frac{w}{G}\Bigr)(\cL w - 
  \vv\cdot \nabla w) \d\xi\ .
$$
Next, using the identity $\cL w = \div(G \nabla (\frac{w}{G}))$ 
and integrating by parts, we obtain
\begin{eqnarray*}
   \inttwo \Bigl(1+\log\frac{w}{G}\Bigr)(\cL w)\d\xi 
   &=& -\inttwo \nabla\Bigl(\log\frac{w}{G}\Bigr)\cdot \frac{G}{w}\nabla 
     \Bigl(\frac{w}{G}\Bigr) w \d\xi \\
   &=& -\inttwo w \Big|\nabla\Bigl(\log\frac{w}{G}\Bigr)\Big|^2 \d\xi 
     \,=\, -I(w)\ .
\end{eqnarray*}
On the other hand, using $\vv\cdot\nabla w = \div (\vv w)$ 
and integrating by parts, we find
\begin{eqnarray*}
  \inttwo \Bigl(1+\log\frac{w}{G}\Bigr)(\vv \cdot \nabla w) \d\xi
    &=& \inttwo (1+\log(4\pi w)) (\vv \cdot \nabla w) \d\xi
    + \inttwo \frac{|\xi|^2}{4}(\vv \cdot \nabla w) \d\xi\ \\
  &=& -\inttwo \vv\cdot\nabla w \d\xi - \half \inttwo (\xxi\cdot\vv)
    w \d\xi .
\end{eqnarray*}
We claim that both integrals in the last expression vanish. This 
is obvious for the first one, since $\vv\cdot\nabla w = \div (\vv w)$.
As for the second one, using \reff{SBS} and Fubini's theorem 
we obtain 
\begin{eqnarray*}
  \inttwo (\xxi\cdot\vv(\xi)) w(\xi) \d\xi &=& 
    \frac{1}{2\pi} \int_{\real^2\times\real^2} \xxi \cdot
    \frac{(\xxi - \eeta)^{\perp}}{|\xi - \eta |^2} w(\eta)
    w(\xi)\d\eta\d\xi \\
  &=& \frac{1}{4\pi} \int_{\real^2\times\real^2} (\xxi-\eeta) \cdot
    \frac{(\xxi - \eeta)^{\perp}}{|\xi - \eta |^2} w(\eta)
    w(\xi)\d\eta\d\xi = 0\ .
\end{eqnarray*}
Summarizing, we have shown that
\begin{equation}\label{Intid}
  H(w(\tau_1)) - H(w(\tau_0)) = - \int_{\tau_0}^{\tau_1}
  I(w(\tau))\d\tau\ ,
\end{equation}
for all $\tau_1 > \tau_0 > 0$. 

We now return to the general case where $w_0 \in L^2(m)\cap
\Sigma_+$. Given $\tau_1 > \tau_0 > 0$, the solution $w$ 
of \reff{SV} satisfies $w \in C^0([\tau_0,\tau_1],H^k(m))$
for any $k \in \intplus$, where $H^k(m)$ is the weighted 
Sobolev space defined in analogy with \reff{H1def}. 
Moreover, the map
$$
   w_0 \in L^2(m) \quad \mapsto \quad w \in 
   C^0([\tau_0,\tau_1],H^k(m))
$$
is continuous. On the other hand, using for instance
\cite{lions:1995}, it is not difficult to verify that the quantity
$I(w)$ is finite for any positive $w \in H^k(m)$ if $k \ge 2$, and
that $I(w)$ depends continuously on $w$ in that topology.  Thus we see
that both sides of \reff{Intid} depend continuously on the initial
data $w_0$ in the topology of $L^2(m)$.  Since \reff{Intid} holds for
all $w_0$ in the dense subset $\cS(\real^2)\cap\Sigma_+$, it follows
that \reff{Intid} is valid for all $w_0 \in L^2(m)\cap\Sigma_+$ and
all $\tau_1 > \tau_0 > 0$.  This concludes the proof. \QED

\subsection{Convergence in $L^2(m)$}\label{convL2m}

Using the compactness properties of Section~\ref{compact_section} and
the two Lyapunov functions of the previous subsection, we are now able 
to prove Proposition~\ref{global}. 

\begin{lemma}\label{a-o} \
Assume that $m > 1$ and that $\{w(\tau)\}_{\tau\in\real}$ is a complete 
trajectory of \reff{SV} which is bounded in $L^2(m)$. Then 
$w(\tau) = \alpha G$ for all $\tau \in \real$, where $\alpha = 
\inttwo w(\xi,0)\d\xi$. 
\end{lemma}

\proof We know from Lemma~\ref{SV_neg} that $\{w(\tau)\}_{\tau\in\real}$
is relatively compact in $L^2(m)$. In view of Remark~\ref{bootstrap}, 
we can assume without loss of generality that $m > 3$. Let $\Omega$
be the $\omega$-limit set of the trajectory $\{w(\tau)\}_{\tau\in\real}$. 
Since by Lemma~\ref{PhiLyap} $\Phi$ is a Lyapunov function which 
is strictly decreasing except on $\Sigma$, it follows from LaSalle's 
invariance principle that $\Omega \subset \Sigma$. In particular, 
since the total mass is conserved, any $\bar w \in \Omega$ satisfies
$\Phi(\bar w) = |\inttwo \bar w(\xi)\d\xi| = |\alpha|$. The same is 
true for any function $\underline{w}$ in the $\alpha$-limit set 
$\cA$. As $\tau \mapsto \Phi(w(\tau))$ is non-increasing, it follows that 
$\Phi(w(\tau)) = |\alpha|$ for all $\tau \in \real$. By Lemma~\ref{PhiLyap}
again, we conclude that $w(\tau) \in \Sigma$ for all $\tau \in \real$. 
Thus, upon replacing $w(\xi_1,\xi_2,\tau)$ by $-w(\xi_2,\xi_1,\tau)$, 
we can assume that $\{w(\tau)\}_{\tau\in\real} \subset 
L^2(m)\cap\Sigma_+$. 

We now use the second Lyapunov function $H$. By Lemma~\ref{HLyap} 
and LaSalle's principle, $\cA$ and $\Omega$ are contained in the 
line of equilibria $\{\alpha' G\}_{\alpha' \ge 0}$. Since the total 
mass is conserved, we necessarily have $\cA = \Omega = \{\alpha G\}$. 
As $H$ is non-increasing, it follows that $H(w(\tau)) = H(\alpha G) =
\alpha \log(\alpha)$ for all $\tau \in \real$. By Lemma~\ref{HLyap}
again, we conclude that $w(\tau) = \alpha G$ for all $\tau \in \real$. 
\QED

\medskip\noindent
{\bf Proof of Proposition~\ref{global}:} 
Let $w_0 \in L^2(m)$ with $m > 1$, and let $w \in
C^0([0,+\infty),L^2(m))$ be the solution of \reff{SV} with initial
data $w_0$. From Lemma~\ref{SV_compact}, we know that
$\{w(\tau)\}_{\tau\ge0}$ is relatively compact in $L^2(m)$. Let
$\Omega \subset L^2(m)$ denote the $\omega$-limit set of this
trajectory. As is well-known, $\Omega$ is non-empty, compact, fully
invariant under the evolution of \reff{SV}, and $\Omega$ attracts
$w(\tau)$ in $L^2(m)$ as $\tau \to +\infty$.  If $\bar w \in \Omega$,
there exists a complete trajectory $\{\bar w(\tau) \}_{\tau\in\real}$
of \reff{SV} such that $\bar w(\tau) \in \Omega$ for all $\tau \in
\real$ and $\bar w(0) = \bar w$. By Lemma~\ref{a-o}, $\bar w(\tau) =
\alpha G$ for all $\tau \in \real$, where $\alpha = \inttwo \bar
w(\xi)\d\xi = \inttwo w_0(\xi)\d\xi$. Thus $\Omega = \{\alpha G\}$,
which is the desired result.
\QED

\subsection{Convergence in $L^1(\real^2)$}\label{convL1}

We now study the behavior of the solutions of \reff{SV} in 
$L^1(\real^2)$ and prove Theorem~\ref{L1glob} and 
Proposition~\ref{Unique}. 

\begin{proposition}\label{L1convprop}
Let $w_0 \in L^1(\real^2)$, and let $w \in C^0([0,\infty),L^1(\real^2))$
be the solution of \reff{SV} with initial data $w_0$. Then 
$|w(\tau)-\alpha G|_1 \to 0$ as $\tau \to \infty$, where 
$\alpha = \inttwo w_0(\xi)\d\xi$. 
\end{proposition}

\proof
We know from Lemma~\ref{L1comp} that $\{w(\tau)\}_{\tau\ge 0}$ 
is relatively compact in $L^1(\real^2)$. Let $\Omega$ be the 
$\omega$-limit set of this trajectory. Then $\Omega$ is non-empty, 
compact, fully invariant under the evolution of \reff{SV}, and $\Omega$ 
attracts $w(\tau)$ in $L^1(\real^2)$ as $\tau \to +\infty$.
If $\bar w \in \Omega$, there exists a sequence $\tau_n 
\to \infty$ such that $|w(\tau_n)-\bar w|_1 \to 0$ and 
$w(\xi,\tau_n) \to \bar w(\xi)$ for almost all $\xi \in \real^2$.
Using \reff{CL1} and Lebesgue's dominated convergence theorem, 
we obtain
\begin{equation}\label{CL4}
  |\bar w(\xi)| = \lim_{n \to \infty} |w(\xi,\tau_n)| 
  \le C_\beta |w_0|_1 \,e^{-\beta|\xi|^2/4}~, \quad 
  \xi \in \real^2~,
\end{equation}
since $a(\tau_n) \to 1$ as $n \to \infty$. In particular, 
this shows that $\Omega$ is bounded in $L^2(m)$ for all 
$m > 1$. 

Now, since $\Omega$ is invariant under the semiflow of \reff{SV},
there exists a complete trajectory $\{\bar w(\tau)\}_{\tau\in \real}$
such that $\bar w(\tau) \in \Omega$ for all $\tau \in \real$ and $\bar
w(0) = \bar w$.  As we just observed, $\{\bar
w(\tau)\}_{\tau\in\real}$ is bounded in $L^2(m)$ for all $m > 1$. 
Applying Lemma~\ref{a-o}, we conclude that $\bar w(\tau) = \alpha
G$ for all $\tau \in \real$, where $\alpha = \inttwo \bar w(\xi)\d\xi
= \inttwo w_0(\xi)\d\xi$.  This proves that $\Omega = \{\alpha G\}$,
which is the desired result.  \QED

\medskip\noindent
{\bf Proof of Theorem~\ref{L1glob}:} 
Let $\omega \in C^0([0,\infty),L^1(\real^2))$ be the solution of
\reff{V} with initial data $\omega_0$. If we set $w(\xi,\tau) = 
e^{\tau}\omega(\xi\,e^{\tau/2},e^{\tau}{-}1)$, then $w \in 
C^0([0,\infty),L^1(\real^2))$ is the solution of \reff{SV} 
with initial data $w_0 = \omega_0$. Applying 
Proposition~\ref{L1convprop} and returning to the original function 
$\omega(x,t)$, we obtain
$$
  \lim_{t \to \infty} \Big|\omega(\cdot,t) -
  \frac{\alpha}{t{+}1}G\Big(\frac{\cdot}{\sqrt{t{+}1}}\Big)
  \Big|_1 = 0~,
$$
which is equivalent to \reff{L1conv} for $p = 1$. Next, interpolating
between \reff{L1conv} for $p = 1$ and \reff{lpvorticity} for 
$p = \infty$, we obtain \reff{L1conv} for $p \in (1,\infty)$.
Then, using \reff{HLS} if $2 < q < \infty$ or Lemma~2.1(b) in 
\cite{gallay:2001a} if $q = \infty$, we arrive at \reff{L1conv2}. 
Finally, using the previous results and the integral equation 
satisfied by $\omega(x,t)$, it is straightforward to show that
\reff{L1conv} also holds for $p = \infty$. \QED

\medskip
As we already observed, a bounded negative trajectory of \reff{SV} in
$L^1(\real^2)$ need not be relatively compact. However, if we assume
that the trajectory is relatively compact or at least that its
$\alpha$-limit set is nonempty, then we have the following result
which generalizes Lemma~\ref{a-o}.

\begin{proposition}\label{L1a-o}
Let $\{w(\tau)\}_{\tau\in\real}$ be a complete trajectory of
\reff{SV} in $L^1(\real^2)$, and assume that $w(\tau)$ has
a convergent subsequence in $L^1(\real^2)$ as $\tau \to 
-\infty$. Then $w(\tau) = \alpha G$ for all $\tau\in \real$, where
$\alpha = \inttwo w(\xi,0)\d\xi$. 
\end{proposition}

\proof
By assumption, there exists $\bar w$ in $L^1(\real^2)$ and 
a sequence $\tau_n \to -\infty$ such that 
$|w(\tau_n)-\bar w|_1 \to 0$ and $w(\xi,\tau_n) \to \bar w(\xi)$ 
for almost all $\xi \in \real^2$. Fix $\tau \in \real$, 
and take $n \ge 0$ sufficiently large so that $\tau_n < \tau$.
In view of \reff{CL1}, we have for all $\xi \in \real^2$
\begin{eqnarray*}
  |w(\xi,\tau)| &\le& C_\beta \inttwo \frac{1}{a(\tau-\tau_n)}\,
    \exp\Bigl(-\beta\frac{|\xi-y e^{-(\tau-\tau_n)/2}|^2}{4
    a(\tau-\tau_n)}\Bigr)|w(y,\tau_n)|\d y \\
  &\le& C_\beta \inttwo \frac{1}{a(\tau-\tau_n)}\,
    \exp\Bigl(-\beta\frac{|\xi-y e^{-(\tau-\tau_n)/2}|^2}{4
    a(\tau-\tau_n)}\Bigr)|w(y,\tau_n) - \bar w(y)|\d y \\
  &+& C_\beta \inttwo \Big| \frac{1}{a(\tau-\tau_n)}\,
    \exp\Bigl(-\beta\frac{|\xi-y e^{-(\tau-\tau_n)/2}|^2}{4
    a(\tau-\tau_n)}\Bigr) - e^{-\beta|\xi|^2/4}\Big| 
    |\bar w(y)|\d y \\
  &+& C_\beta |\bar w|_1 e^{-\beta|\xi|^2/4}~.
\end{eqnarray*}
Taking the limit $n \to \infty$ and using Lebesgue's dominated
convergence theorem, we obtain
$$
  |w(\xi,\tau)| \le C_\beta |\bar w|_1 e^{-\beta|\xi|^2/4}~,
  \quad \xi \in \real^2\,,~\tau \in \real~.
$$
This shows that the trajectory $\{w(\tau)\}_{\tau \in \real}$
is bounded in $L^2(m)$ for any $m > 1$, and the result follows
from Lemma~\ref{a-o}. \QED

\medskip
As is clear from the change of variables \reff{omega-w}, 
results about negative trajectories of \reff{SV} give 
information on the behavior of solutions of \reff{V} 
as $t \to 0+$. In particular, Lemma~\ref{a-o} and 
Proposition~\ref{L1a-o} show that solutions of \reff{V}
with Dirac masses as initial data are unique in a certain class. 
A generalization of these results allows to prove 
Proposition~\ref{Unique}. 

\medskip\noindent
{\bf Proof of Proposition~\ref{Unique}:} 
In view of Theorem~\ref{L1V2}, we can assume without loss of
generality that $T = \infty$. The solution $\omega(x,t)$ of 
\reff{V} can be represented as
$$
  \omega(x,t) = \inttwo \Gamma_u(x,t;y,s)\omega(y,s)\d y~,
  \quad x \in \real^2~, \quad t > s > 0~.
$$
Here $\Gamma_u$ is the fundamental solution of the time-dependent
linear operator $\partial_t - \Delta + \uu\cdot\nabla$, and 
$\uu(x,t)$ is the velocity field obtained from $\omega(x,t)$ 
via the Biot-Savart law. By assumption, there exists $K > 0$ such 
that $|\omega(\cdot,t)|_1 \le K$ for all $t > 0$. From 
\cite{carlen:1996}, we know that for any $\beta \in (0,1)$
there exists $C_\beta > 0$ (depending on $K$) such that 
\begin{equation}\label{gamm1}
  |\Gamma_u(x,t;y,s)| \le \frac{C_\beta}{t-s}
  \,\exp\Bigl(-\beta\frac{|x-y|^2}{4(t-s)}\Bigr)\ ,
\end{equation}
for all $x,y \in \real^2$ and all $t > s > 0$, see \reff{pointbdd}. 
Moreover, it is shown in \cite{osada:1987} (see also Theorem~3.1 
in \cite{giga:1988}) that $\Gamma_u$ is a H\"older continuous
function of its arguments. More precisely, there exists 
$\gamma \in (0,1)$ (depending only on $K$) and, for any $\tau > 0$, 
a constant $C > 0$ (depending only on $K$ and $\tau$) such that 
that
\begin{equation}\label{gamm2}
  |\Gamma_u(x,t;y,s) - \Gamma_u(x,t;y',s')| \le C
  \Bigl(|y-y'|^\gamma + |s-s'|^{\gamma/2}\Bigr)\ ,
\end{equation}
whenever $t-s\ge\tau$ and $t-s'\ge\tau$. In particular, if $x,y
\in\real^2$ and $t > 0$, the function $s \mapsto \Gamma_u(x,t;y,s)$
can be continuously extended to $s = 0$, and this extension 
(still denoted by $\Gamma_u$) satisfies \reff{gamm1}, \reff{gamm2}
with $s = 0$. 

Now, fix $x \in \real^2$ and $t > 0$. Then for any $s \in (0,t)$
we have
\begin{eqnarray*}
  \omega(x,t) &=& \inttwo \Gamma_u(x,t;y,0)\omega(y,s)\d y \\
  &+& \inttwo (\Gamma_u(x,t;y,s) - \Gamma_u(x,t;y,0))\omega(y,s)\d y~.
\end{eqnarray*}
In view of \reff{gamm2}, the second integral in the right-hand side
converges to zero as $s \to 0+$, since $|\omega(\cdot,s)|_1 \le K$
for all $s > 0$. On the other hand, since $y \mapsto 
\Gamma_u(x,t;y,0)$ is continuous and bounded, and since 
$\omega(\cdot,s)$ converges weakly to $\alpha\delta$ as $s \to 0+$, 
the first integral converges to $\alpha \Gamma_u(x,t;0,0)$. Thus
$$
  |\omega(x,t)| = |\alpha| |\Gamma_u(x,t;0,0)| \le \frac{C_\beta
  |\alpha|}{t}\,e^{-\beta|x|^2/(4t)}\ ,\quad x \in \real^2\, 
  \quad t > 0\ .
$$
Finally, let $w(\xi,\tau) = e^\tau\omega(\xi e^{\tau/2},e^\tau)$ 
for $\xi \in \real^2$, $\tau \in \real$. Then $w \in C^0(\real,
L^1(\real^2))$ is a solution of \reff{SV} which satisfies 
$|w(\xi,\tau)| \le C_\beta |\alpha|\,e^{-\beta|\xi|^2/4}$
for all $\xi \in \real^2$, $\tau \in \real$. In particular, 
$\{w(\cdot,\tau)\}_{\tau\in\real}$ is bounded in $L^2(m)$ 
for any $m > 1$, hence by Lemma~\ref{a-o} $w(\xi,\tau) = 
\alpha' G(\xi)$ for some $\alpha' \in \real$. Clearly 
$\alpha' = \alpha$, and the proof is complete. \QED

\begin{remark}\label{unique2}
As is clear from the proof, Proposition~\ref{Unique} remains
true if one assumes only that $\omega(\cdot,t)$ stays bounded 
in $L^1(\real^2)$ and that $\omega(\cdot,t_n)$ converges 
weakly to $\alpha \delta$ for some sequence $t_n \to 0$.  
\end{remark}

\begin{remark}\label{localstab}
A slight extension of the techniques developped in this section
allows to prove that the Oseen vortices are stable in the
sense of Lyapunov\: for any $\epsilon > 0$, there exists
$\delta > 0$ such that, for all initial data $w_0 \in 
L^1(\real^2)$ with $|w_0 - \alpha G|_1 \le \delta$, the
solution of \reff{SV} satisfies $|w(\tau) - \alpha G|_1 \le 
\epsilon$ for all $\tau \ge 0$. Note that this does not 
follow from Proposition~\ref{L1convprop}. 
\end{remark}

\subsection{Convergence rate for positive solutions}

If we restrict ourselves to {\em nonnegative} solutions of
\reff{SV}, then combining the entropy dissipation law 
\reff{H-dot} with a few classical inequalities in information
theory we can obtain an explicit estimate on the time needed
for a solution to approach the Oseen vortex. This is the 
so-called ``entropy dissipation method'', which is by now
a classical approach in kinetic theory, see for instance
\cite{toscani:1999,arnold:2001,villani:2002}. 

Let $w \in L^1(\real^2)$, $w \ge 0$, and assume that 
$\alpha = \inttwo w(\xi)\d\xi > 0$. In information theory, 
the quantity $H(w)$ defined in \reff{Hdef} is called the
{\em relative Kullback entropy} of $w$ with respect to
the Gaussian $G$. Similarly, $I(w)$ defined in \reff{Idef}
is called the {\em relative Fisher information} of $w$ with 
respect to $G$. The difference between the entropy $H(w)$ 
and its minimal value $H(\alpha G)$ is bounded from 
below by the Csisz\'ar-Kullback inequality
\begin{equation}\label{CKineq}
 \frac{1}{2\alpha} \|w - \alpha G\|_{L^1}^2 
 \le H(w) - H(\alpha G)\ ,  
\end{equation}
and from above by the Stam-Gross logarithmic Sobolev 
inequality
\begin{equation}\label{logSob}
  H(w) - H(\alpha G) \le I(w)\ .
\end{equation}

Assume now that $w_0 \in L^2(m) \cap \Sigma_+$ for some 
$m > 3$, and that $\alpha = \inttwo w_0(\xi)\d\xi > 0$. 
Let $w \in C^0([0,\infty),L^2(m))$ be the solution of \reff{SV} 
with initial data $w_0$. Combining \reff{H-dot} with \reff{logSob}, 
we immediately obtain
$$
   H(w(\tau)) - H(\alpha G) \le  \bigl(H(w_0) - H(\alpha G)\bigr)
   \,e^{-\tau}\ , \quad \tau \ge 0\ .
$$
Applying \reff{CKineq}, we conclude that
\begin{equation}\label{explicitbd}
  \|w(\tau) - \alpha G\|_{L^1} \le \sqrt{2\alpha}\, 
  \bigl(H(w_0) - H(\alpha G)\bigr)^{1/2}\,e^{-\tau/2}\ ,  
  \quad \tau \ge 0\ .
\end{equation}
This shows that $w(\tau)$ converges to $\alpha G$ at the rate 
$e^{-\tau/2}$, which is optimal in general (see 
Section~\ref{local_section}). Moreover, \reff{explicitbd} 
gives an explicit estimate of the prefactor in terms 
of the initial data. In particular, this provides an 
explicit upper bound of the time needed for the solution 
to enter a given neighborhood of the vortex.

Unfortunately, we do not know how to extend the entropy 
dissipation method to the general case where the vorticity
may change sign. In the next section, we obtain local 
convergence rates by studying \reff{SV} in a neighborhood of 
the family of Oseen vortices, but this approach does
not provide any explicit estimate in the sense of \reff{explicitbd}.


\section{Local convergence rates}\label{local_section}

From the results of the previous section we know that
any solution of the Navier-Stokes equation whose initial vorticity 
distribution lies in $L^1(\real^2)$ will converge toward the Oseen 
vortex with the same total vorticity. In the present section we 
show that for solutions in the weighted space $L^2(m)$ with
$m > 1$ we can derive precise estimates on the rate at which solutions 
approach the vortices. 

Our analysis proceeds by first analyzing the linearization of 
\reff{SV} at a vortex solution. We prove estimates on the location 
of the spectrum of linearized operator which in particular 
imply that the vortex solutions are spectrally stable for all
values of the circulation Reynolds number. We then show that these
bounds also imply decay estimates for the full nonlinear evolution 
in a neighborhood of the vortex. 

\subsection{Eigenvalue estimates}\label{spectrum}

Fixing $\alpha \in \real$ and linearizing \reff{SV} around $w= \alpha G$, 
we find
\begin{equation}\label{linear_eq}
  \partial_t w + \alpha(\vv^G \cdot \nabla w + \vv \cdot\nabla G)
  \,=\, \Delta w + {1 \over 2}\xi\cdot\nabla w + w~,
\end{equation} 
where as usual $\vv$ is the velocity field associated to
$w$ via \reff{SBS} and $\vv^G$ is the velocity field
of the Oseen vortex. 
This equation can be rewritten as $\partial_t w = \cL w - \alpha 
\Lambda w$, where
\begin{equation}\label{linear_ops}
   \cL w \,=\, \Delta w + {1 \over 2}\xi\cdot\nabla w + 
   w~, \quad \Lambda w = \vv^G \cdot \nabla w + 
   \vv\cdot\nabla G~,
\end{equation}
see \reff{R-evol}, \reff{LamRdef}. 

The linear operator $\cL$ in $L^2(m)$ is studied in detail in 
(\cite{gallay:2001a}, Appendix~A). It is defined on the 
maximal domain
$$
  \cD_m(\cL) = \Big\{w \in L^2(m)\,\Big|\, \Delta w + {1 \over 2}
  \xi\cdot\nabla w \in L^2(m)\Big\}\ .
$$
If $w \in \cD_m(\cL)$, one can show that $\Delta w \in L^2(m)$, 
so that $\cD_m(\cL) \subset H^2(m)$. The essential spectrum 
of $\cL$ is given by
$$
   \sigma_m^\ess(\cL) = \Bigl\{\lambda \in \complex \,\Big|\, 
   \Re(\lambda) \le \frac{1-m}{2} \Bigr\}\ .
$$
In addition, $\cL$ has a sequence of eigenvalues $0,\,-1/2,\,-1,\dots$
whose eigenfunctions are rapidly decreasing at infinity.  

Because of the spatial decay of $\vv^G$ and $G$, the operator
$\Lambda$ is a relatively compact perturbation of $\cL$ and hence
$\sigma_m^\ess(\cL-\alpha\Lambda)$ does not depend on $\alpha$.  In
particular, we can always push this essential spectrum far away from
the imaginary axis by taking $m > 0$ sufficiently large. Thus the
spectral stability of the vortex solutions will be determined by the
isolated eigenvalues of $\cL - \alpha \Lambda$ in $L^2(m)$.  As we
shall see, the corresponding eigenfunctions have a Gaussian decay at
infinity so that, in contrast to the essential spectrum, these
isolated eigenvalues do not depend on $m$.

We next observe that, due to symmetries of equation \reff{V},
some eigenvalues of $\cL-\alpha\Lambda$ are in fact independent 
of $\alpha$. For instance, if $m > 1$, then $\lambda = 0$ is a simple 
eigenvalue of $\cL$ in $L^2(m)$, with eigenfunction $G$. Since 
$\vv^G \cdot\nabla G = 0$, it is clear from \reff{linear_ops} that 
$\Lambda G = 0$, so that $0$ is an eigenvalue of $\cL - \alpha 
\Lambda$ for any $\alpha \in \real$. This zero eigenvalue is due to 
the fact that the Oseen vortices form a one-parameter family 
of equilibria of \reff{SV}. The associated spectral projection 
$P_0$ reads (for any $\alpha$) 
$$
   (P_0 w)(\xi) = G(\xi) \inttwo w(\xi')\d\xi'\ .
$$
Thus, it will be sufficient to study the spectrum of $\cL - \alpha 
\Lambda$ in the spectral subspace $L^2_0(m)$ defined in
\reff{L20mdef}, which by \reff{masscons} is also invariant under 
the nonlinear evolution \reff{R-evol}. 

Similarly, if $m > 2$, $\lambda = -1/2$ is a double eigenvalue of
$\cL$ with eigenfunctions $F_1, F_2$, where $F_j = -\partial_j G$
($j = 1,2$). Differentiating the identity $\vv^G \cdot\nabla G = 0$
with respect to $\xi_j$, we see that $\Lambda F_j = 0$ ($j = 1,2$).
It follows that $-1/2$ is still an eigenvalue of $\cL - \alpha 
\Lambda$ for any $\alpha \in \real$. This eigenvalue originates in 
the translation invariance of \reff{V} with respect to $\xi \in 
\real^2$. The associated spectral projection $P_1$ reads 
(for any $\alpha$) 
$$
   (P_1 w)(\xi) = F_1(\xi) \inttwo \xi_1' w(\xi')\d\xi' + 
                  F_2(\xi) \inttwo \xi_2'w(\xi')\d\xi'\ . 
$$
Thus, it is again sufficient to study the spectrum of $\cL - \alpha 
\Lambda$ in the spectral subspace $L^2_1(m)$ defined by
\begin{equation}\label{L21mdef}
  L^2_1(m) \,=\, \Bigl\{w \in L^2_0(m) \,\Big|\, 
  \inttwo \xi_j w(\xi)\d\xi = 0 \hbox{ for } j = 1,2\Bigr\}\ ,
\end{equation}
which by \reff{momcons} is also invariant under the nonlinear
evolution \reff{R-evol}.

Finally, if $m > 3$, $\lambda = -1$ is a triple eigenvalue of
$\cL$ with eigenfunctions $\Delta G$, $(\partial_1^2 -\partial_2^2)G$, 
and $\partial_1 \partial_2 G$. Since $\Delta G = \frac{1}{4} 
(|\xi|^2-4)G$ is radially symmetric, it is clear that 
$\Lambda(\Delta G) = 0$, so that $-1$ is still an eigenvalue of 
$\cL - \alpha \Lambda$ for any $\alpha \in \real$. This is
due the fact that \reff{V} is autonomous and invariant under the 
rescaling $\omega(x,t) \mapsto \lambda^2 \omega(\lambda x, 
\lambda^2 t)$. However, as we shall see, the eigenvalue $-1$ is 
{\em simple} if $\alpha \neq 0$. The associated spectral 
projection $P_2$ reads (for any $\alpha \neq 0$) 
$$
   (P_2 w)(\xi) = \Delta G(\xi) \inttwo \frac{1}{4}
   (|\xi'|^2-4) w(\xi')\d\xi'\ .
$$
Thus, if $\alpha \neq 0$, it is sufficient to study the spectrum 
of $\cL - \alpha \Lambda$ in the spectral subspace $L^2_2(m)$ 
defined by
\begin{equation}\label{L22mdef}
  L^2_2(m) \,=\, \Bigl\{w \in L^2_1(m) \,\Big|\, 
  \inttwo |\xi|^2 w(\xi)\d\xi = 0\Bigr\}\ ,
\end{equation}
which (as can be verified by a direct calculation) is also invariant 
under the nonlinear evolution \reff{R-evol}.

\medskip 
The principal result of this subsection is:

\begin{proposition}\label{eigenvalue_bound}
Fix $m > 1$ and $\alpha \in \real$. Then any eigenvalue $\lambda$  
of $\cL - \alpha\Lambda$ in $L^2_0(m)$ satisfies
\begin{equation}\label{eigg1}
  \Re(\lambda) \le \max\Bigl(-\half\,,\,\frac{1-m}{2}\Bigr)\ .
\end{equation}
If moreover $m > 2$, then any eigenvalue $\lambda$  
of $\cL - \alpha\Lambda$ in $L^2_1(m)$ satisfies
\begin{equation}\label{eigg2}
  \Re(\lambda) \le \max\Bigl(-1\,,\,\frac{1-m}{2}\Bigr)\ .
\end{equation}
\end{proposition}

\begin{remark}\label{sharp} In view of the preceding remarks
estimates \reff{eigg1} and \reff{eigg2} are sharp. If 
$m > 3$ and $\alpha \neq 0$, the proof shows that any eigenvalue
of $\cL - \alpha\Lambda$ in $L^2_2(m)$ satisfies $\Re(\lambda) < -1$, 
but we are not able to give a sharp estimate in that case. 
Numerical calculations in \cite{prochazka:1995} indicate that
the eigenvalues that are not frozen by symmetries have a 
real part that converges to $-\infty$ as $|\alpha| \to \infty$, 
thereby suggesting that a fast rotation has a stabilizing effect 
on the vortex. Proposition~\ref{eigenvalue_bound} shows rigorously 
that at the spectral level perturbations of the vortex solutions 
decay at least as fast when $\alpha$ is large as when $\alpha=0$.
\end{remark}

To prove Proposition~\ref{eigenvalue_bound}, we proceed in 
three steps. First, we observe that the linear operators 
$\cL$ and $\Lambda$ are invariant under the group of 
rotations $SO(2)$. Thus, using polar coordinates in 
$\real^2$ and expanding the angular variable in Fourier
series, we reduce the eigenvalue equation for the operator
$\cL - \alpha\Lambda$ to a (nonlocal) ordinary differential 
equation in the radial variable. Next, a careful study of
this differential equation reveals that, if $\lambda \in 
\complex$ is an isolated eigenvalue of $\cL - \alpha\Lambda$, 
the corresponding eigenfunction has a Gaussian decay at
infinity. Finally, we prove that the operators $\cL$ and
$\Lambda$ are respectively self-adjoint and skew-symmetric
in a weighted $L^2$ space with appropriate Gaussian 
weight, and Proposition~\ref{eigenvalue_bound} then follows 
from elementary considerations. 

\subsubsection{Polar coordinates}\label{polar}

Fix $m > 0$. For any $n \in \allint$, let $P_n$ be the orthogonal 
projection in $L^2(m)$ defined by
\begin{eqnarray}\nonumber
 &&(P_n w)(r\cos\theta,r\sin\theta) = \omega_n(r) e^{in\theta}~,
  \quad {\rm where}\\ \label{omegandef}
 &&\omega_n(r)= \frac{1}{ 2\pi}\int_0^{2\pi} w
  (r\cos\theta,r\sin\theta) e^{-in\theta}~\d\theta~.
\end{eqnarray}
Clearly, $P_n P_{n'} = \delta_{nn'}P_n$ and $\sum_{n \in \allint}P_n 
= \one$. If $w \in L^2(m)$, the functions $\omega_n : 
\real_+ \to \complex$ in \reff{omegandef} belong to the Hilbert space
\begin{equation}\label{Zmdef}
   Z(m) \,=\, \Bigl\{\omega : \real_+ \to \complex \,\Big|\, 
   \int_0^\infty r(1+r^2)^m |\omega(r)|^2 \d r < \infty\Bigr\}~.
\end{equation}
For any $n \in \allint$, let $\cL_n$ be the linear operator on $Z(m)$ 
defined by
\begin{equation}\label{cLndef}
   \cL_n \omega \,=\, \partial_r^2 \omega + \Bigl({r \over 2}+{1 \over r}
   \Bigr)\partial_r \omega + \Bigl(1 - {n^2 \over r^2}\Bigr) \omega~.
\end{equation}
Let also $\Lambda_n$ be the (bounded) linear operator on $Z(m)$ defined
by $\Lambda_0 = 0$ and
\begin{equation}\label{Lambdandef}
   \Lambda_n \omega \,=\, in(\phi \omega - g\Omega)~, \quad n \neq 0~,
\end{equation}
where
$$
   \phi(r) \,=\, {1 \over 2\pi r^2}(1-e^{-r^2/4})~, \quad
   g(r) \,=\, {1 \over 4\pi}\,e^{-r^2/4}~,
$$
and
\begin{equation}\label{Omegadef}
   \Omega(r) \,=\, {1 \over 4|n|}\left(\int_0^r \Bigl({z \over r}
   \Bigr)^{|n|} z\omega(z)\d z + \int_r^\infty \Bigl({r \over z}
   \Bigr)^{|n|} z\omega(z)\d z\right)~.
\end{equation}
It is easy to see that the operator $\Lambda_n$ is indeed well-defined:

\begin{lemma}\label{welldef}
If $n \in \allint^*$ and $\omega \in Z(m)$ for some $m > 0$, 
then \reff{Omegadef} defines a continuous function $\Omega : 
\real_+ \to \complex$. Moreover, $\Omega(r)/r$ converges to 
zero as $r \to \infty$ and is at most logarithmically divergent
as $r \to 0$. 
\end{lemma}

\proof The proof is straightforward using \reff{Zmdef}, \reff{Omegadef}
and H\"older's inequality.\QED

\medskip
We now show that $\cL_n$ and $\Lambda_n$ are the expressions of 
$\cL$ and $\Lambda$ in polar coordinates: 

\begin{lemma}\label{commute_op}
The operators $\cL$ and $\Lambda$ commute with the projections 
$P_n$. If $n \in \allint$ and $w \in \cD_m(\cL)$ for some $m > 0$, 
then
\begin{eqnarray}\label{PnL}
   (\cL P_n w)(r\cos\theta,r\sin\theta) &=&
   e^{in\theta} (\cL_n \omega_n)(r)~, \\ \label{PnLambda}
   (\Lambda P_nw)(r\cos\theta,r\sin\theta) &=& 
   e^{in\theta} (\Lambda_n \omega_n)(r)~,
\end{eqnarray}
where $\omega_n(r) = e^{-in\theta} (P_nw)(r\cos\theta,r\sin\theta)$.
\end{lemma}

\proof All we need is to prove \reff{PnL} and \reff{PnLambda}. 
The first relation follows from \reff{linear_ops} by an elementary
calculation. To prove \reff{PnLambda}, assume that 
$w(r\cos\theta,r\sin\theta) = \omega_n(r)e^{in\theta}$ for some 
$n \in \allint$ and some $\omega_n \in Z(m)$. Then
$$
   \vv^G\cdot\nabla w \,=\, {1 \over 2\pi r}(1-e^{-r^2/4}) 
   {1 \over r}\partial_\theta w \,=\, ine^{in\theta} 
   \phi \omega_n~.
$$
On the other hand, the velocity field $v$ corresponding to $w$ 
satisfies $\partial_1 v_2 -\partial_2 v_1= w$, $\partial_1 v_1 + 
\partial_2 v_2 = 0$. In polar coordinates, these relations become
$$
   {1 \over r}\partial_r(rv_\theta) -{1 \over r}\partial_\theta v_r 
   \,=\, w~, \quad {1 \over r}\partial_r(rv_r) + {1 \over r}
   \partial_\theta v_\theta \,=\, 0~.
$$
We look for a solution of the form $v_r = \bar v_r(r)\,e^{in\theta}$,
$v_\theta = \bar v_\theta(r)\,e^{in\theta}$. Then
$$
   (r\bar v_\theta)' -in \bar v_r \,=\, r\omega_n~, \quad
   (r\bar v_r)' + in\bar v_\theta \,=\, 0~.
$$
Eliminating $\bar v_\theta$, we find the following ODE for 
$h = r\bar v_r$:
$$
   (rh')' -{n^2 \over r}h + in r \omega_n \,=\, 0~.
$$
The general solution is:
$$
   h(r) \,=\, {in \over 2|n|}\left(\int_0^r \Bigl({z \over r}
   \Bigr)^{|n|} z\omega_n(z)\d z + \int_r^\infty \Bigl({r \over z}
   \Bigr)^{|n|} z\omega_n(z)\d z\right) + A_1 r^n + A_2 r^{-n}~,
$$
where $A_1, A_2 \in \complex$. Since we want a velocity 
$\bar v_r = h/r$ that is locally integrable and 
converges to zero at infinity, we must choose $A_1=A_2=0$. Setting 
$\Omega = {1 \over 2in}h = {1 \over 2in} r\bar v_r$, we finally 
obtain:
$$
   \vv\cdot\nabla G \,=\, -{1 \over 2}rv_r g \,=\, -in e^{in\theta}
   g\Omega~.
$$
This concludes the proof of \reff{PnLambda}. \QED

\subsubsection{Gaussian decay of eigenfunctions}

The aim of this paragraph is to prove:

\begin{lemma}\label{eigen_decay} Fix $m > 0$, and assume that $w \in 
\cD_m(\cL)$ satisfies $(\cL -\alpha\Lambda)w = \mu w$, where 
$\alpha \in \real$ and $\Re(\mu) > {1-m \over 2}$. Then there exist 
$C > 0$ and $\gamma \ge 0$ such that
$$
   |w(\xi)| \le C(1+|\xi|^2)^\gamma \,e^{-|\xi|^2/4}\ ,\quad 
   \xi \in \real^2\ .
$$
\end{lemma}

\proof We use the notations of the preceding paragraph. 
According to Lemma~\ref{commute_op}, we can assume that the 
eigenfunction $w$ satisfies $w = P_n w$ for some $n \in \allint$.
Thus, there exists $\omega \in Z(m)$ such that $w(r\cos\theta,
r\sin\theta) = \omega(r) e^{i n \theta}$. In view of \reff{cLndef}, 
\reff{Lambdandef} and Lemma~\ref{commute_op}, $\omega$ satisfies
the (inhomogeneous) ordinary differential equation
\begin{equation}\label{radialone}
\omega''(r) + \left( \frac{r}{2} +\frac{1}{r} \right) \omega'(r)
+ \left( 1 - \mu -\frac{n^2}{r^2} - i n \alpha \phi \right) \omega(r)
+ i n \alpha g \Omega(r) = 0 \ ,
\end{equation}
where $\Omega$ is defined in \reff{Omegadef}. 

The basic idea is now to use the classical results of Coddington
and Levinson \cite{coddington:1955} to show that, for
$r$ large, any solution of \reff{radialone} either decays
like $\omega(r) \sim r^{2 \gamma} e^{-r^2/4}$ for some
$\gamma \ge 0$, or like $\omega(r) \sim 
r^{2\mu-2}$ in which case it cannot belong to $L^2(m)$. 
However, a certain amount of preliminary work is needed in order 
to bring \reff{radialone} into a form to which we can apply the 
results of \cite{coddington:1955}.

We begin by introducing new variables $f$, $F$ and $t$ via the
definitions
$$
  \omega(r) = f(r^2/4)\ , \ 
  \Omega(r) = F(r^2/4)\ ,\ t=r^2/4\ .
$$
In terms of these new variables, \reff{radialone} takes the form
\begin{equation}\label{radial}
  f''(t) + \left(1+\frac{1}{t}\right)f'(t)
  +\left(\frac{1-\mu}{t} - a(t)\right) f(t) + b(t) = 0
  \ ,\ \  t> 0\ ,
\end{equation}
where
$$
  a(t) = \frac{n^2}{4 t^2} + \frac{i n \alpha}{8 \pi t^2} (1-e^{-t})
  \ ,\quad  b(t) = \frac{i n \alpha}{4 \pi t} e^{-t} F(t)\ .
$$

Recall that we are interested in the behavior of solutions of this
equation for $t$ large.  We first consider the behavior of the
homogeneous part of this equation and then construct the solution
of the full equation via the method of variation of parameters.

\begin{lemma}\label{homogeneous_asym} The linear, homogeneous
equation 
$$
  \tf''(t) + \left(1+\frac{1}{t}\right)\tf'(t)
  +\left(\frac{1-\mu}{t} - a(t)\right) \tf(t) = 0\ ,
$$
has two linearly independent solutions $\phi_1(t)$ and
$\phi_2(t)$ such that 
$$
  \lim_{t\to \infty} t^{1-\mu} \left( { {\phi_1(t)}\atop{{\phi_1}'(t)}
  }\right) \,=\, \left( {1}\atop{0} \right)\ ,\quad \lim_{t\to \infty} 
  t^{\mu}e^t \left( { {\phi_2(t)}\atop{{\phi_2}'(t)}} \right)
  \,=\, \left( {1}\atop{-1} \right)\ .
$$
\end{lemma}

\proof If we define $x(t) = \left( { {\tf(t)}\atop{{\tf}'(t)}} \right)$
we can rewrite the differential equation in the lemma as
$$
  x'(t) = (A +V(t) +R(t)) x(t)\ ,
$$
where
$$
  A= \pmatrix{0&1 \cr 0&-1\cr}\ ,\ 
  V(t) = \pmatrix{0&0 \cr -\frac{1-\mu}{t} & -\frac{1}{t} }\ ,\ 
  R(t) = \pmatrix{0&0 \cr a(t) &0 \cr}\ .
$$
But in this form, the lemma follows immediately from 
(\cite{coddington:1955}, Theorem 3.8.1, p. 92).
\QED

\medskip
We now derive the asymptotic form of the solution $f(t)$
of \reff{radial} by applying the method of variation
of parameters.  Set 
\begin{equation}\label{varf}
f(t) = A(t) \phi_1(t) + B(t) \phi_2(t)\ ,
\end{equation}
where $A'(t) \phi_1(t) +B'(t) \phi_2(t) =0$.  Then differentiating
$f$ and using \reff{radial} we find
$$
  A'(t) \phi_1'(t) + B'(t) \phi_2'(t) + b(t) = 0\ .
$$
Solving for $A'$ and $B'$ we obtain
$$
  \left( {A'(t)}\atop{B'(t)} \right) =
  \frac{1}{W(t)} \left( {b(t) \phi_2(t)}\atop{-b(t) \phi_1(t)} \right)\ ,
$$
where $W(t) = -\frac{1}{t} e^{-t}$ is the Wronskian determinant 
of $\phi_1$ and $\phi_2$. Integrating
both sides of this equation we find
$$
  A(t) = A_1 - \int_1^t s e^s b(s) \phi_2(s) \d s\ , \quad 
  B(t) = B_1 + \int_1^t se^s b(s) \phi_1(s) \d s\ .
$$
Recall from the definition of $b(t)$ that 
$e^s b(s) = \frac{i n \alpha}{4 \pi} \frac{F(s)}{s}$.  
From Lemma~\ref{welldef} we know that $F(s)/s = \Omega(2\sqrt{s})/s$ 
converges to zero as $s \to \infty$, so that $e^s b(s)$ is bounded 
for $s \ge 1$. Thus, the asymptotic
behavior of $\phi_2$ implies that $A_1(t) \to A_{\infty}$
as $t \to \infty$.  If $A_{\infty} \ne 0$, then $f(t) \sim
A_{\infty} t^{\mu-1}$ as $t \to \infty$ and hence (reverting
to the original polar coordinates) $\omega(r) \sim
\omega_{\infty} r^{2\mu -2}$.  But since $\Re(\mu) > \frac{1-m}{2}$,
this would imply that 
$\int_0^{\infty} r(1+r^2)^m |\omega(r)|^2 \d r = \infty$ and
this in turn would violate the hypothesis that the eigenfunction
$w$ is in $L^2(m)$.  Thus, $A_{\infty} =0$, and 
$A(t) = \int_t^{\infty} s e^s b(s) \phi_2(s) \d s$ from which
we conclude that $|A(t)| \le C e^{-t} t^{\gamma}$, for
some $\gamma \ge 0$.  In analogous fashion one proves that
$|B(t)| \le C t^{\gamma}$.  Inserting these bounds on $A$
and $B$ into \reff{varf} and using the asymptotic estimates
on $\phi_1$ and $\phi_2$ we conclude that there exists
$\gamma \ge 0$ such that
$|f(t)| \le C t^{\gamma} e^{-t}$, for $t \ge 1$, or
$$
   |\omega(r)| \le C r^{2 \gamma} \,e^{-r^2/4}\ ,\ \ r \ge 1\ .
$$
This is the desired estimate, since $|w(\xi)| = |\omega(|\xi|)|$.
The proof of Lemma~\ref{eigen_decay} is now complete. \QED

\subsubsection{Localization of eigenvalues}

Let $X$ denote the (complex) Hilbert space
$$
   X = \Bigl\{w\in L^2(\real^2) \,\Big|\, G^{-1/2}w \in 
   L^2(\real^2)\Bigr\}~,
$$
equipped with the scalar product
$$
  (w_1,w_2)_X \,=\, \inttwo {1 \over G(\xi)}
  \,\bar w_1(\xi) w_2(\xi) \d \xi~.
$$
We also introduce the closed subspaces $X_0, X_1$ defined by
\begin{eqnarray*}
   X_0 &=& \Big\{w \in X \,\Big|\, \inttwo w(\xi)\d \xi \,=\, 0\Big\}
    \,=\, X \cap L^2_0(m)\ ,\\
   X_1 &=& \Big\{w \in X_0 \,\Big|\, \inttwo \xi_j w(\xi)\d \xi \,=\, 0
    \hbox{ for }j = 1,2\Big\} \,=\, X \cap L^2_1(m)\ .
\end{eqnarray*}
According to Lemma~\ref{eigen_decay}, if $\mu$ is an 
eigenvalue of $\cL -\alpha\Lambda$ in $L^2(m)$ with 
$\Re(\mu) > \frac{1-m}{2}$, the corresponding eigenfunction
belongs to $X$. This result is very useful because both operators 
$\cL$ and $\Lambda$ have nice properties in this space.

\begin{lemma}\label{Lsymm} The linear operator $\cL$ is self-adjoint 
in $X$, and $\cL \le 0$. Morevoer, $\cL \le -1/2$ on $X_0$ and 
$\cL \le -1$ on $X_1$. 
\end{lemma}

\proof Define $L : \cD(L) \to L^2(\real^2)$ by $\cD(L) = 
\{\psi \in H^2(\real^2)\,|\, |\xi|^2\psi \in L^2(\real^2)\}$ and
$$
   L \,=\, G^{-1/2}\cL \,G^{1/2} \,=\, \Delta -{|\xi|^2 \over 16} 
   + {1 \over 2}~.
$$
In quantum mechanics, the operator $-L$ is (up to numerical constants)
the Hamiltonian of the harmonic oscillator in $\real^2$. As is 
well-known (see for example \cite{glimm:1989}), $L$ is self-adjoint 
in $L^2(\real^2)$ and $\sigma(L) = \{-n/2\,|\,n \in \intplus\}$. 
By construction, the operator $\cL : \cD(\cL) \to X$ with domain
$$
  \cD(\cL) = \Big\{w \in X\,\Big|\, |\xi|^2 w \in X\ ,~\Delta w 
  + {1\over 2}\xi\cdot\nabla w \in X\Big\}\ ,
$$
is thus self-adjoint in $X$ with the same spectrum. In particular, 
$\cL \le 0$. Now, observe that $0$ is a simple eigenvalue of $\cL$
with eigenfunction $G$, and that $X_0$ is just the orthogonal
complement of the eigenspace $\real G$ in $X$. Thus $X_0$ is stable
under $\cL$ and the restriction of $\cL$ to $X_0$ is a self-adjoint
operator satisfying $\cL \le -1/2$. Similarly, one can show that $\cL
\le -1$ on $X_1$. \QED

\begin{lemma}\label{Lamasymm} The linear operator $\Lambda$ is 
skew-symmetric in $X$.
\end{lemma}

\proof
Since $\Lambda G = 0$ and since the subspace $X_0$ is stable 
under $\Lambda$, it is sufficient to show that $\Lambda$ is 
skew-symmetric on $X_0$. Let $w, \tilde w \in X_0 \cap 
\cD(\cL)$, and denote by $\vv, \tilde \vv$ the corresponding 
velocity fields. Without loss of generality, we assume that 
$w, \tilde w$ (hence also $\vv, \tilde \vv$) are real functions. 
Then 
$$
   (\tilde w,\Lambda w)_X = \inttwo
   \left({1 \over G}\,\tilde w \vv^G\cdot\nabla w -{1 \over 2}
   \tilde w (\vv \cdot \xi)\right)\d \xi~,
$$
because $\nabla G = -{\xi \over 2}G$. Observe that $G^{-1}\vv^G$ is 
a divergence free vector field, so that
\begin{equation}\label{eqA}
   \inttwo {1 \over G}\,\tilde w \vv^G \cdot \nabla w\d \xi = 
   -\inttwo {1 \over G}\,w \vv^G \cdot \nabla \tilde w \d \xi~.
\end{equation}
On the other hand, the following identity is easy to check:
$$
   \tilde w (\vv \cdot \xi) + w (\tilde \vv \cdot \xi) \,=\,
   (\xi_1\partial_1 -\xi_2\partial_2)(v_1 \tilde v_2 + v_2 \tilde v_1)
   + (\xi_1 \partial_2 + \xi_2 \partial_1) (v_2 \tilde v_2 - 
   v_1 \tilde v_1)~.
$$
Since $w, \tilde w$ have zero mean, it follows from 
(\cite{gallay:2001a}, Appendix~B) that $\vv, \tilde \vv$ 
decay at least like $1/|\xi|^2$ as $\xi \to \infty$. Thus, integrating 
both sides, we obtain
\begin{equation}\label{eqB}
   \inttwo (\tilde w (\vv\cdot \xi) + w (\tilde \vv \cdot \xi))
   \d \xi \,=\, 0~.
\end{equation}
Combining \reff{eqA} and \reff{eqB}, we see that 
$(\tilde w,\Lambda w)_X + (\Lambda \tilde w, w)_X = 0$. \QED

\medskip\noindent
{\bf Proof of Proposition~\ref{eigenvalue_bound}.}
Fix $m > 1$ and assume that $\lambda$ is an eigenvalue of
$\cL-\alpha\Lambda$ in $L^2_0(m)$ with $\Re(\lambda) > \frac{1-m}{2}$. 
By Lemma~\ref{eigen_decay}, there exists a nonzero $w \in X_0
\cap \cD(\cL)$ such that $(\cL-\alpha\Lambda)w = \lambda w$. 
In particular, 
$$
   \lambda(w,w)_X \,=\, (w,\cL w)_X -\alpha(w,\Lambda w)_X~,
$$
hence
$$
   \Re(\lambda) (w,w)_X \,=\, (w,\cL w)_X \,\le\, -{1 \over 2}
   (w,w)_X~,
$$
since $\Lambda$ is skew-symmetric and $\cL \le -1/2$ on $X_0$. 
Thus, $\Re(\lambda) \le -1/2$. 

Similarly, if $m > 2$ and $\lambda$ is an eigenvalue of
$\cL-\alpha\Lambda$ in $L^2_1(m)$ with $\Re(\lambda) > \frac{1-m}{2}$,
there exists a nonzero $w \in X_1 \cap \cD(\cL)$ such that 
$(\cL-\alpha\Lambda)w = \lambda w$. Proceeding as above and 
using the fact that $\cL \le -1$ on $X_1$, we obtain 
$\Re(\lambda) \le -1$. \QED

\begin{remark}\label{more_eig}
If $m > 3$ and $\lambda$ is an eigenvalue of $\cL-\alpha\Lambda$ 
in $L^2_2(m)$ with $\Re(\lambda) > \frac{1-m}{2}$, there exists a 
nonzero $w \in \cD(\cL) \cap L^2_2(m)$ such that $(\cL-\alpha\Lambda)w 
= \lambda w$. The argument above shows that $\Re(\lambda) \le -1$, 
and that $\Re(\lambda) = -1$ if and only if $\cL w = -w$. But
this implies that $w$ is a linear combination of $(\partial_1^2 
-\partial_2^2)G$ and $\partial_1 \partial_2 G$, and a direct 
calculation shows that no such $w$ can be an eigenfunction of 
$\cL-\alpha\Lambda$ if $\alpha \neq 0$. Thus any eigenfunction 
$\lambda$ of $\cL-\alpha\Lambda$ in $L^2_2(m)$ satisfies 
$\Re(\lambda) < -1$ if $\alpha \neq 0$. 
\end{remark}

\subsection{Bounds on the linear evolution}
\label{linear_subsection}

In this subsection we prove that the eigenvalue estimates
of the previous subsection imply analogous bounds on
the linear evolution. Fix $\alpha \in \real$ and consider
the linear equation $\partial_\tau R = \cL R -\alpha \Lambda R$
which is the linearization of \reff{SV} a the vortex $\alpha G$,
see \reff{R-evol}, \reff{linear_eq}. The corresponding integral
equation reads
\begin{equation}\label{eqintR}
  R(\tau) = \cS(\tau) R_0 - \alpha
 \int_0^{\tau} \cS(\tau-s) \Lambda R(s) \d s\ ,
\end{equation}
where $\cS(\tau) = \exp(\tau\cL)$. Proceeding as in 
Section~\ref{compactL2}, it is straightforward to
show by a contraction mapping argument that this equation defines 
a strongly continuous semigroup $\cT_\alpha(\tau)$ in $L^2(m)$ 
for any $m > 1$, namely $R(\tau) = \cT_\alpha(\tau)R_0$. We first 
prove that $\cT_\alpha(\tau)$ is a compact perturbation of $\cS(\tau)$. 

\begin{lemma}\label{K_compact} Let $m > 1$. The linear operator 
$\cK_\alpha(\tau) = \cT_\alpha(\tau) - \cS(\tau)$ is compact in $L^2(m)$ 
for each $\tau > 0$.
\end{lemma}

\proof
All the necessary estimates are already contained in 
Section~\ref{compactL2}. Observe that the term 
$\cK_\alpha(\tau)R_0$ in \reff{eqintR} is precisely what we called
$R_2(\tau) + R_3(\tau)$ in the proof of Lemma~\ref{SV_compact}.
Repeating the estimates proved there, we obtain
$$
   \|\cK_\alpha(\tau)R_0\|_{m+1} \le C(\tau)\|R_0\|_m\ ,\quad \tau > 0\ ,
$$
for some $C(\tau) > 0$. Similarly, the result of Lemma~\ref{SV_smooth}
applies to $\cT_\alpha(\tau)R_0$, hence to $\cK_\alpha(\tau)R_0$, 
and yields
$$
   \|\nabla \cK_\alpha(\tau)R_0\|_{m} \le \frac{C(\tau)}{a(\tau)^{1/2}}
   \,\|R_0\|_m\ ,\quad \tau > 0\ ,
$$
where $a(\tau) = 1-e^{-\tau}$. The conclusion then follows from 
Rellich's criterion. \QED

\medskip
We now use this result to bound the essential spectrum of
$\cT_\alpha(\tau)$. We recall that $\lambda \in \complex$ is 
in the essential spectrum of a linear operator $L$ if $\lambda$
is not a {\em normal point} for $L$, i.e. if $\lambda$ is not 
in the resolvent set of $L$ and is not an isolated eigenvalue 
of $L$ of finite multiplicity, see e.g. \cite{henry:1981}.
For any $r > 0$, we denote by $\cB(r)$ the closed disk of 
radius $r$ centered at the origin in $\complex$:  
$$
   \cB(r) \,=\, \{z \in \complex \,|\, |z| \le r\}\ .
$$

\begin{lemma}\label{T_ess} Let $m > 1$. For any $\tau > 0$, 
the essential spectrum of $\cT_\alpha(\tau)$ in $L^2(m)$ satisfies
\begin{equation}\label{sigma_ess}
   \sigma_m^\ess(\cT_\alpha(\tau)) \,=\, \cB(e^{\tau(1-m)/2})\ .
\end{equation}
\end{lemma}

\proof
Fix $\tau > 0$. The results of (\cite{gallay:2001a}, Appendix~A) imply
that the essential spectrum of $\cS(\tau)$ in $L^2(m)$ is exactly
$\cB(e^{\tau(1-m)/2})$. Since $\cT_\alpha(\tau)$ is a compact perturbation of
$\cS(\tau)$, it follows from (\cite{henry:1981}, Theorem A.1) that the
complement of the $\cB(e^{\tau(1-m)/2})$ in $\complex$ either consists
entirely of eigenvalues of $\cT_\alpha(\tau)$, or entirely of normal points
for $\cT_\alpha(\tau)$. In our case, the first possibility is excluded.
Indeed, assume that $\lambda \in \complex$ is an eigenvalue of $\cT_\alpha(\tau)$
with $|\lambda| > 1$. By the spectral mapping theorem, $\lambda = 
e^\nu$ where $\nu$ is an eigenvalue of $\cL -\alpha\Lambda$ in
$L^2(m)$ with $\Re(\nu) > 0$, which contradicts
Proposition~\ref{eigenvalue_bound}. Thus, $\sigma_m^\ess(\cT_\alpha(\tau))
\subset \cB(e^{\tau(1-m)/2})$, and since $\cS(\tau)$ is also a compact
perturbation of $\cT_\alpha(\tau)$ the same argument shows that
$\sigma_m^\ess(\cT_\alpha(\tau)) = \cB(e^{\tau(1-m)/2})$.
\QED

\medskip
By construction, the spectral subspaces $L^2_0(m)$ and 
$L^2_1(m)$ are left invariant by the semigroup $\cT_\alpha(\tau)$. 
Combining Lemma~\ref{T_ess} with the eigenvalue estimates 
of Section~\ref{spectrum}, we obtain precise bounds on the 
growth of $\cT_\alpha(\tau)$ in these subspaces: 

\begin{proposition}\label{decay} Fix $\alpha \in \real$. 
Assume that either\\
{\bf a)} $0 < \mu \le 1/2$ and $R_0 \in L^2_0(m)$ for some 
  $m > 1+2\mu$, or\\
{\bf b)} $1/2 < \mu \le 1$ and $R_0 \in L^2_1(m)$ for some 
  $m > 1+2\mu$.\\
Then there exists $C > 0$ (independent of $R_0$) such that
\begin{equation}\label{decay_est}
  \|\cT_\alpha(\tau)R_0\|_m \le C \,e^{-\mu \tau} \|R_0\|_m\ , \quad 
  \tau \ge 0\ .
\end{equation}
\end{proposition}

\proof Assume first that $0 < \mu < 1/2$. By Lemma~\ref{T_ess},
the essential spectrum of $\cT_\alpha(1)$ in $L^2_0(m)$ satisfies 
$\sigma_m^\ess(\cT_\alpha(1)) = \cB(e^{(1-m)/2})$, and $e^{(1-m)/2} <
e^{-\mu}$ since $m > 1+2\mu$. If $\lambda \in \sigma(\cT_\alpha(1))$
satisfies $|\lambda| > e^{(1-m)/2}$, then $\lambda$ is 
an isolated eigenvalue of $\cT_\alpha(1)$ and (by the spectral 
mapping theorem) there exists an eigenvalue $\nu$ of $\cL -\alpha 
\Lambda$ in $L^2_0(m)$ such that $e^\nu = \lambda$. Applying 
Proposition~\ref{eigenvalue_bound}, we obtain $\Re(\nu) \le 1/2$, 
hence $|\lambda| \le e^{-1/2} < e^{-\mu}$. Thus the spectral 
radius of $\cT_\alpha(1)$ in $L^2_0(m)$ is strictly less than 
$e^{-\mu}$, and \reff{decay_est} follows (see e.g. 
\cite{engel:2000}, Proposition IV.2.2). A similar argument 
proves \reff{decay_est} if $R_0 \in L^2_1(m)$ and $1/2 < \mu < 1$.

Now, assume that $\mu = 1/2$ and $m > 2$. Any $R_0 \in L^2_0(m)$
can be decomposed as $R_0 = \beta_1 F_1 + \beta_2 F_2 + \tilde R_0$, 
where $\beta_j = \inttwo \xi_j R_0 \d\xi$, $F_j = -\partial_j G$, 
and $\tilde R_0 \in L^2_1(m)$. It follows that
$$
  \cT_\alpha(\tau) R_0 = e^{-\tau/2}(\beta_1 F_1 + \beta_2 F_2)
  + \cT_\alpha(\tau) \tilde R_0\ ,
$$
and we already know that $\|\cT_\alpha(\tau) \tilde R_0\|_m 
\le C e^{-\nu\tau} \|\tilde R_0\|_m$ for some $\nu > 1/2$. 
Thus \reff{decay_est} holds and is sharp in this case. 
A similar argument shows that \reff{decay_est} holds with 
$\mu = 1$ if $R_0 \in L^2_1(m)$ for some $m > 3$. \QED

\medskip
Finally, we will need in the following subsection $L^p$-$L^q$ 
estimates of $\cT_\alpha(\tau)$ and its derivatives, in the spirit
of \reff{LpLq}. 

\begin{proposition}\label{A5-analogue} Under the assumptions
of Proposition~\ref{decay}, if $R_0$ satisfies in addition 
$b^m R_0 \in L^q(\real^2)$ for some $q \in (1,2)$, then
$$
  \|\cT_\alpha(\tau) R_0\|_m \le \frac{C\,e^{-\mu\tau}}{a(\tau)^{\frac{1}{q}
  -\frac{1}{2}}} \,|b^m R_0|_q \ , \quad
  \|\nabla \cT_\alpha(\tau) R_0\|_m \le \frac{C\,e^{-\mu\tau}}{a(\tau)^{\frac{1}{q}}}
  \,|b^m R_0|_q \ , \quad \tau > 0\ .
$$
\end{proposition}

\proof
Fix $1 < q < 2$ and assume that $b^m R_0 \in L^q(\real^2)$. Given 
$T > 0$, we consider the function space
\begin{eqnarray*}
  Y &=& \{R \in C^0((0,T],H^1(m)) \,|\, \|R\|_Y < \infty\}\ , \quad
    \hbox{where }\\
  \|R\|_Y &=& \sup_{0 < \tau \le T} \Bigl\{a(\tau)^{\frac{1}{q}
   -\frac{1}{2}}\|R(\tau)\|_m + a(\tau)^{\frac{1}{q}}
   \|\nabla R(\tau)\|_m\Bigr\}\ .
\end{eqnarray*}
Using \reff{LpLq} and proceeding as in the proof of Lemma~\ref{SV_smooth}, 
it is straightforward to show that, if $T > 0$ is sufficiently 
small, there exists $C > 0$ such that \reff{eqintR} has a unique 
solution in the ball of radius $C|b^m R_0|_q$ centered at the 
origin in $Y$. This proves the desired estimates for $0 < \tau 
\le T$, and the general case follows from Proposition~\ref{decay}
if one uses the semigroup property and the fact that the subspaces 
$L^2_0(m)$ and $L^2_1(m)$ are left invariant by $\cT_\alpha(\tau)$. 
\QED

\subsection{Bounds on the nonlinear evolution}
\label{nonlinear_subsection}

In this subsection, we show that our estimates on the linear 
semigroup $\cT_\alpha(\tau)$ generated by $\cL -\alpha\Lambda$
imply similar bounds on the full nonlinear evolution \reff{R-evol}. 
As is easy to verify, the integral equation \reff{Rint} satisfied
by $R(\tau)$ is equivalent to
\begin{equation}\label{Rint2}
  R(\tau) = \cT_\alpha(\tau) R_0 - \int_0^{\tau} \cT_\alpha(\tau -s) 
  (\vv^{R}(s) \cdot \nabla R(s))\d s\ .
\end{equation}
The following result implies Proposition~\ref{local1} as a 
particular case: 

\begin{proposition}\label{local2}
Under the assumptions of Proposition~\ref{decay}, the solution
$R(\tau)$ of \reff{R-evol} with initial data $R_0$ satisfies
$\|R(\tau)\|_m = \cO(e^{-\mu\tau})$ as $\tau \to \infty$. 
\end{proposition}

\proof
Let $R \in C^0([0,+\infty),L^2_0(m))$ be the solution of 
\reff{R-evol} with initial data $R_0$. We know from 
Proposition~\ref{global} that $\|R(\tau)\|_m$ converges to zero 
as $\tau \to \infty$, hence we can assume without loss of 
generality that $\|R_0\|_m$
is small. Given $T > 0$, we define
$$
  M(T) = \max\Bigl\{\,\sup_{0 \le \tau \le T} \,e^{\mu\tau}\|R(\tau)\|_m 
   \,,\, \sup_{0 < \tau \le T} a(\tau)^{\half} \,e^{\mu\tau}
   \|\nabla R(\tau)\|_m\Bigr\}\ .
$$
By Lemma~\ref{SV_smooth}, $M(T) < \infty$ and there exists 
$C_0 > 0$ such that $M(T) \le C_0 \|R_0\|_m$ if $T$ is 
sufficiently small. Applying Propositions~\ref{decay} and 
\ref{A5-analogue}, and using the fact that the subspaces 
$L^2_0(m)$ and $L^2_1(m)$ are stable under the nonlinearity 
$\vv^R \cdot \nabla R$, we obtain from \reff{Rint2}
$$
   \|R(\tau)\|_m \le C \,e^{-\mu\tau} \|R_0\|_m\ + C \int_0^{\tau} 
   e^{-\mu(\tau-s)} \frac{1}{a(\tau{-}s)^{\frac{1}{q}-\frac{1}{2}}}
   |b^m \vv^{R}(s)\cdot\nabla R(s)|_q \d s\ ,
$$
where $1 < q < 2$. As in the proof of Lemma~\ref{SV_smooth}, we 
can bound
$$
   |b^m \vv^{R}(s)\cdot\nabla R(s)|_q \le C \|R(s)\|_m \|\nabla R(s)\|_m
   \le \frac{C\,e^{-2\mu s}}{a(s)^\frac{1}{2}} \,M(T)^2\ ,
$$
hence
\begin{eqnarray*}
   \|R(\tau)\|_m &\le& e^{-\mu\tau} \Bigl(C \|R_0\|_m\ + C M(T)^2 
   \int_0^{\tau} \frac{e^{-\mu s}}{a(\tau{-}s)^{\frac{1}{q}-\frac{1}{2}}
     a(s)^\frac{1}{2}} \d s\Bigr) \\
   &\le& e^{-\mu\tau} \Bigl(C_1 \|R_0\|_m\ + C_2 M(T)^2\Bigr)\ , \quad 
   0 \le \tau \le T\ ,
\end{eqnarray*}
where $C_1, C_2 > 0$ are independent of $T$. Without loss of 
generality, we assume in what follows that $C_1 \ge C_0$. 
Differentiating \reff{Rint2} and using similar estimates, we also obtain
$$
   \|\nabla R(\tau)\|_m \le \frac{e^{-\mu\tau}}{a(\tau)^\frac{1}{2}} 
   \Bigl(C_1 \|R_0\|_m\ + C_2 M(T)^2\Bigr)\ , \quad 0 < \tau \le T\ .
$$
Summarizing, we have shown 
\begin{equation}\label{trinome}
  M(T) \le C_1 \|R_0\|_m\ + C_2 M(T)^2\ . 
\end{equation}
Now, assume that $R_0$ is small enough so that $4 C_1 C_2 \|R_0\|_m < 1$.
As $C_1 \ge C_0$, we also have $C_0 \|R_0\|_m \le \bar M$, 
where
$$
   \bar M = \frac{1}{2C_2}\Bigl(1 - \sqrt{1-4 C_1 C_2 \|R_0\|_m}\Bigr)
   > 0\ .
$$
Thus $M(T) \le \bar M$ for $T > 0$ sufficiently small, and since $M(T)$
depends continuously on $T$ it follows from \reff{trinome} that 
$M(T) \le \bar M$ for all $T > 0$. In particular, $\|R(\tau)\|_m \le 
\bar M\,e^{-\mu\tau}$ for all $\tau \ge 0$. \QED

\begin{remark}\label{faster_conv}
Similarly, if $R_0 \in L^2_2(m)$ for some $m > 3$, the solution 
$R(\tau)$ of \reff{R-evol} with initial data $R_0$ satisfies
$\|R(\tau)\|_m = \cO(e^{-\mu\tau})$ as $\tau \to \infty$, for
some $\mu \ge 1$ depending on $\alpha$. We know from 
Remark~\ref{more_eig} that $\mu > 1$ if $\alpha \neq 0$, but 
we have no sharp result in that case. 
\end{remark}

\begin{remark}\label{localstab2}
It follows directly from the proof that the Oseen vortices 
are stable equilibria in $L^2(m)$, see Remark~\ref{localstab}.
\end{remark}

We conclude this section by showing that, in the case where the total
vorticity is nonzero, there is no loss of generality in assuming that
the perturbations of the vortex have vanishing first order
moments. Indeed, assume that $1/2 < \mu \le 1$ and that $w_0 \in
L^2(m)$ for some $m > 1+2\mu$. Let $w \in C^0([0,\infty),L^2(m)$ be
the solution of \reff{SV} with initial data $w_0$. For any $b \in
\real^2$, the function $\bar w(\xi,\tau) = w(\xi +
b\,e^{-\tau/2},\tau)$ is again a solution of \reff{SV} (because the
original equation \reff{V} is translation invariant in $x \in
\real^2$). If $\alpha = \inttwo w_0(\xi)\d\xi \neq 0$, we can choose
$(b_1,b_2) = \alpha^{-1}(\beta_1,\beta_2)$, where $\beta_i = \inttwo
\xi_i w_0(\xi)\d\xi$. Then $\inttwo \xi_i \bar w(\xi,0)\d\xi = 0$ for
$i = 1,2$, so that $\bar w(\cdot,0) - \alpha G \in L^2_1(m)$. Applying
Proposition~\ref{local2}, we obtain $\|\bar w(\tau)-\alpha G\|_m =
\cO(e^{-\mu\tau})$ as $\tau \to \infty$. Returning to the original
function $w(\xi,\tau)$ and using a straightforward Taylor expansion,
we obtain the second order asymptotics
\begin{equation}\label{second}
   \|w(\tau) - \alpha G -(\beta_1 F_1 + \beta_2 F_2)\,e^{-\tau/2}\|_m
   = \cO(e^{-\mu\tau}), \quad \tau \to \infty\ ,
\end{equation}
where $F_j = -\partial_j G$, $j = 1,2$. As was already mentioned, 
this result still holds when $\alpha = 0$ except that, if $\mu = 1$
and $(\beta_1,\beta_2) \neq (0,0)$, the right-hand side of 
\reff{second} should be replaced by $\cO(\tau e^{-\tau})$. 

\begin{remark} Assume that $m > 3$, $\alpha \neq 0$, and 
$\inttwo |\xi|^2 w_0(\xi)\d\xi \neq 0$. If we translate 
and rescale $w_0$ appropriately, we can produce a new initial
condition $\bar w_0$ with $\bar w_0 - \alpha G \in L^2_2(m)$.
By Remark~\ref{faster_conv}, the corresponding solution then 
satisfies $\|\bar w(\tau)-\alpha G\|_m = \cO(e^{-\mu\tau})$ as 
$\tau \to \infty$ for some $\mu > 1$. Moreover $w(\xi,\tau)$ and
$\bar w(\xi,\tau)$ are linked by a simple relation, due to
the fact that the original equation \reff{V} is translation and 
dilation invariant. Using this relation, we obtain that the  
the next correction to the asymptotics \reff{second} is of the 
form $\gamma (\Delta G) \,e^{-\tau}$, for some $\gamma \in \real$. 
\end{remark}

\begin{remark}
The connection of the translation invariance of the Navier-Stokes
equation with the decay associated to the first moment of the
vorticity seems first to have been remarked upon by Bernoff and
Lingevitch, \cite{bernoff:1994}. The connection between symmetries of
the linear and nonlinear heat equation and Burgers' equation and the
decay rates of the long-time asymptotics of solutions of these
equations was systematically explored in \cite{witelski:1998} and
\cite{miller:2003}.
\end{remark}

\bibliographystyle{plain}

\end{document}